\documentclass[11pt]{article}
\title{On analogues of B\" acklund theorem \\ in affine differential geometry of surfaces}
\author{Maria Robaszewska\\
\textit{\small{Instytut Matematyki, Uniwersytet Pedagogiczny,}}\\
\textit{\small{ul. Podchor\c a\.zych 2, 30-084 Krak\' ow, Poland.}}\\
\textit{\small{e-mail: maria.robaszewska@up.krakow.pl}}}
\date{ \ }

\usepackage{amsmath}
\usepackage{amssymb}

\newcommand*{\R}{\mathbf{R}}
\newcommand*{\Proj}{\mathbf{P}}
\newcommand*{\im}{\mathrm{im}}
\newcommand*{\linsp}{\mathrm{span}}

\newcommand*{\Ric}{\mathrm{Ric}}
\newcommand*{\sign}{\mathrm{sign}}

\newcommand*{\n}{\mathbf{n}}

\newtheorem{theorem}{Theorem}[section]
\newtheorem{lemma}[theorem]{Lemma}
\newtheorem{proposition}[theorem]{Proposition}

\begin{document}
\maketitle

{\small
\textbf{Abstract:} 
We recall the well-known Chern--Terng theorem concerning affine minimal surfaces. Next we formulate some complementary (with transversal fields necessarily not parallel) affine B\" acklund theorem. We describe some geometrical conditions which imply the local symmetry of both induced connections. We give also some necessary and sufficient conditions under which the affine fundamental forms are proportional.
\medskip

\medskip

}

\section{Introduction.}
The most classical B\" acklund theorem is the following B\" acklund theorem for surfaces in Euclidean space:

\begin{theorem}
Let $f,\widehat f:M\to \R ^3$, be a pair of surfaces in the Euclidean space $\R ^3$
satisfying the following conditions:
\medskip

\noindent \emph{(i)} for every $p\in M$ $f(p)\ne \widehat f(p)$,
the vector $\widehat f(p)-f(p)$ is tangent to $f(M)$ at $f(p)$ and is tangent to $\widehat f(M)$ at $\widehat f(p)$,
\medskip

\noindent \emph{(ii)} the length $L:=|\widehat f(p)-f(p)|$ of $\widehat f(p)-f(p)$ is independent of $p$,
\medskip

\noindent \emph{(iii)} the angle $\sigma $ between the Euclidean normals $\n $ and $\widehat \n$ (of $f(M)$ and $\widehat f(M)$ respectively) is constant and $\sin \sigma \ne 0$.
\medskip

Then both surfaces are of constant negative Gaussian curvature $\kappa =\widehat \kappa =-\frac{\sin ^2\sigma }{L^2}$.
The second fundamental forms $h$ and $\widehat h$ of $f$ and $\widehat f$ are proportional.
\end{theorem}

In this article we will present analogues of B\" acklund theorem in affine differential geometry of surfaces. 
We recall Chern--Terng theorem and prove some other affine B\" acklund theorem, concerning surfaces with locally symmetric induced connection. 

Our aim was to generalize B\" acklund theorem to the situation, when in ambient space there is only the volume form, and we cannot measure length or angle. The volume form is parallel with respect to the standard linear connection $D$ in $\R ^3$.
We study two immersions $f$ and $\widehat f$, which are focal surfaces of a rectilinear congruence. Each of them is endowed with an equiaffine transversal vector field, $\xi $ and $\widehat {\xi }$ respectively. 
Unlike the Euclidean normals, those transversal fields are not determined by the immersions. Of course, one may use the affine normal, and this particular case will be also considered. 
We will impose on $(f,\xi )$ and $(\widehat f,\widehat {\xi })$ some conditions which guarantee that both induced connections $\nabla $ and $\widehat {\nabla }$ are locally symmetric. Our idea was to consider the volume of the parallelepiped spanned by $\widehat f-f$ and both transversal fields. In Euclidean case this volume is a non-zero constant. 
The conjecture that condition of constant volume together with some other conditions about the values of conormal map enforce both Blaschke connections to be locally symmetric turned out to be true. Some partial result, with Blaschke normal of $f$ tangent to $\widehat f$ and vice versa, is contained in 
\cite{MR}. 
However, in case of arbitrary equiaffine transversal fields $\xi $ and $\widehat {\xi }$ one should admit also non-constant volume $\det (\widehat f-f, \xi ,\widehat {\xi }) $.

Our result seems to be a common generalization of the classical B\" acklund theorem (see for example \cite{ChT} or \cite{T}) and Minkowski space B\" acklund theorem (\cite{AB}, \cite{MR}). It also includes the case of non-metrizable connections with $\dim \im R=1$, studied by Opozda in \cite{O}. The theorem is complementary to Chern and Terng analogue of B\" acklund's theorem in affine geometry \cite{ChT}, because in 
\cite{ChT} the affine normals $\R \xi $ and $\R \widehat {\xi }$ were assumed to be parallel, hence $\det (\widehat f-f, \xi ,\widehat {\xi }) =0$.

\section{Preliminaries}

We recall the basic notions of affine differential geometry. More details can be found in \cite{NS}. Here we consider only two-dimensional manifolds immersed into affine space $\R ^3$. The standard connection in $\R ^3$ is denoted by $D$.
\smallskip

Let $f:M\to \R ^3$ be an immersion of a two-dimensional manifold $M$ into $\R ^3$. Let $\xi :M\to \R ^3$ be a transversal vector field. For each $p\in M$ we have the decomposition $\R ^3=f_*(T_pM)\oplus \R \xi _p$. The \emph{induced connection} $\nabla $, the \emph{affine fundamental form} $h$ (relative to the transversal vector field $\xi $), the \emph{affine shape operator} $S$ and the \emph{transversal connection form} $\tau $ are defined by the following Gauss and Weingarten formulae
\begin{equation}
D_Xf_*(Y)=f_*(\nabla _XY)+h(X,Y)\, \xi ,
\end{equation}
\begin{equation}
D_X\xi =-\, f_*(SX)+\tau (X)\, \xi .
\end{equation}
The volume element induced by $(f,\xi )$ on $M$ is
\[ \theta (X,Y)=\det (f_*X,f_*Y,\xi ).\]
The determinant $\det _{\theta }h$ of a symmetric covariant tensor $h$ of degree $2$ relative to $\theta $ is, by definition, equal to $\det [h_{ij}]$, where $h_{ij}=h(X_i,X_j)$ and $X_1$, $X_2$ is a unimodular basis for $\theta $: $\theta (X_1,X_2)=1$.
Let $(\R ^3)^*$ be the dual space of the vector space $\R ^3$. For immersion $f:M\to \R ^3$ with transversal vector field $\xi :M\to \R ^3$ the \emph{conormal map} $\nu :M\to (\R ^3)^*$ is defined as follows:
\[ \nu _p(f_*(X_p)):=0 \ \ \mbox{and} \ \ \nu _p(\xi _p):=1 \ \ \mbox{for} \ \ p\in M, \ \ X_p\in T_pM.\]
\smallskip

The rank of the affine fundamental form is independent of the choice of transversal vector field. If $h$ is nondegenerate, then we say that the surface is \emph{nondegenerate}.
If $f$ is nondegenerate, then for each point $p\in M$ there exists a transversal vector field defined in a neighbourhood of $p$ satisfying the conditions

(I) $\nabla \theta =0$

(II) $\theta $ coincides with the volume element of the nondegenerate metric $h$.

Such a transversal vector field is unique up to a sign and is called the \emph{affine normal field} or \emph{Blaschke normal field}. The connection induced by the affine normal vector field is called the \emph{Blaschke connection} and $h$ is called the \emph{affine metric}. The condition (I) is equivalent to $\tau =0$ and the condition (II) is equivalent to $|\det _{\theta }h|=1$. An {\it equiaffine} transversal field is a transversal field satisfying the condition $\tau =0$.
\bigskip

B\" acklund theorem is usually formulated for two focal surfaces of some rectilinear congruence. A rectilinear congruence is a two-parametric family of straight lines. Under some additional assumption about the congruence one can find two families of ruled developable surfaces with rulings belonging to the congruence. Each line of the congruence is contained in one developable surface of each family and is tangent to the edge of regression of this developable surface at the point which is called the {\it focal point}. Except of some particular degenerate cases the set of all focal points forms two {\it focal surfaces}. We parametrize the focal surfaces in such a way that $f(p)$ and $\widehat f(p)$ belong to the same straight line of congruence. We may consider the mapping $f(p)\mapsto \widehat f(p)$ between the two focal surfaces. If this mapping preserves the asymptotic lines, a rectilinear congruence is called a $W$-{\it congruence}.
\smallskip

More details about rectilinear congruences one can find for example in \cite{E}.

\section{A necessary and sufficient condition \\ for rectilinear congruence with non-degenerate focal surfaces to be a W-congruence}

In this section we will study the condition that the affine fundamental forms $h$ and $\widehat h$, of $(f,\xi )$
and $(\widehat f,\widehat {\xi })$ respectively, are proportional. In Euclidean or Minkowski space B\" acklund theorem this condition is a part of the assertion, whereas in affine case it is an assumption. 
\medskip

\begin{proposition}
\label{funkcja}
Let $f:M\to \R ^3$ and $\widehat f:M\to \R ^3$ be
non-degenerate immersions of a two-dimensional real manifold $M$ into
affine space $\R ^3$ such that for every $p\in M$ $f(p)\ne \widehat f(p)$,
the vector $\widehat f(p)-f(p)$ is tangent to $f(M)$ at $f(p)$ and is tangent to $\widehat f(M)$ at $\widehat f(p)$.
\smallskip

Let $\xi $ and $\widehat {\xi }$ be some transversal vector fields for $f$ and $\widehat f$ respectively. We denote by $h$ and $\widehat h$ the corresponding affine fundamental forms, and by $\nu $ and $\widehat {\nu }$ the conormal maps.
\medskip

Then:
\smallskip

\noindent \emph{(i)} If $\det (\widehat f-f, \xi , \widehat {\xi })=0$, \ then \ $1-\nu (\widehat {\xi })\, \widehat {\nu }(\xi )=0$.
\medskip

\noindent \emph{(ii)} If $1-\nu (\widehat {\xi })\, \widehat {\nu }(\xi )=0$ at some point $p$ and $f_*(T_pM)\ne \widehat f_*(T_pM)$, then $\det (\widehat f(p)-f(p), \xi _p, \widehat {\xi _p})=0$.

\noindent \emph{(iii)} If $\xi $ and $\widehat {\xi }$ are such that $\det (\widehat f-f,\xi ,\widehat {\xi })\ne 0$ and $\widetilde {\xi }=\lambda \, \xi +f_*Z$, \linebreak $\overline {\xi }=\mu \, \widehat {\xi }+\widehat f_*V$, then
\[ 1-\widetilde {\nu }(\overline {\xi })\, \overline {\nu }(\widetilde {\xi })=
\frac{1-\nu (\widehat {\xi })\, \widehat {\nu }(\xi )}{\lambda \, \mu \, \det (\widehat f-f,\xi ,\widehat {\xi })}\, \det (\widehat f-f,\widetilde {\xi },\overline {\xi }).\]

\noindent \emph{(iv)} If moreover $\det (\widehat f-f,\widetilde {\xi },\overline {\xi })\ne 0$, then
\[ \left ( \frac{1-\widetilde {\nu }(\overline {\xi })\,
\overline {\nu }(\widetilde {\xi })}{\det (\widehat f-f,\widetilde {\xi },\overline {\xi })}\right ) ^4\frac{1}{{\det }_{\widetilde {\theta }}\widetilde h\ {\det }_{\overline {\theta }}\overline h}=\left ( \frac{1-\nu (\widehat {\xi })\,
\widehat {\nu }(\xi )}{\det (\widehat f-f,\xi ,\widehat {\xi })}\right ) ^4\frac{1}{{\det }_{\theta }h\ {\det }_{\widehat {\theta }}\widehat h}.\]
\end{proposition}

\emph{Proof.} (i)
There exist nowhere vanishing vector fields $X_1$ and $\widehat X_1$ on $M$ such that
\begin{equation}
\widehat f-f=f_*X_1
\label{fgx1}
\end{equation}
and
\begin{equation}
\widehat f-f=\widehat f_*\widehat X_1.
\label{fgx1dasz}
\end{equation}
Since $f_*X_1$ and $\xi $ are linearly independent, from $\det (\widehat f-f,\xi ,\widehat {\xi })=0$ it follows that $\widehat {\xi }=\alpha \, f_*X_1+\beta \, \xi $ for some $\alpha $ and $\beta $. Here $\beta \ne 0$, because $f_*X_1=\widehat f_*\widehat X_1$ is tangent to $\widehat f$.
We have $\nu (\widehat {\xi })=\beta $ and from $\widehat {\xi }=\alpha \, \widehat f_*\widehat X_1+\beta \, \xi $ we obtain $1=\beta \, \widehat {\nu }(\xi )$.
\smallskip

(ii) Conversely, if $\nu (\widehat {\xi })\, \widehat {\nu }(\xi )=1$, then $\widehat {\xi }=f_*T+A\, \xi $ and $\xi =\widehat f_*U+\frac{1}{A}\, \widehat {\xi }$ with some $A\ne 0$. It follows that $f_*T=\widehat {\xi }-A\, \xi =-\, \widehat f_*(A\, U)$. Therefore $f_*T$ is tangent to $f$ and is tangent to $\widehat f$. By assumption $f_*T_pM\ne \widehat f_*T_pM$, hence $f_*T_pM\cap \widehat f_*T_pM=\R f_*{X_1}_p$ and $\widehat {\xi }\in \linsp \{ f_*{X_1}_p,\xi _p\} $.
\smallskip

(iii) Let $W=\det (\widehat f-f,\xi ,\widehat {\xi })$, $A=\nu (\widehat {\xi })$ and $\widehat A=\widehat {\nu }(\xi )$.
For every $p\in M$, $\dim f_*T_pM=2$, $\dim \linsp \{ \xi _p,\widehat {\xi }_p\} =2$ and $f_*T_pM\ne \linsp \{ \xi _p,\widehat {\xi }_p\} $, because $\xi _p\notin f_*T_pM$. Therefore $\dim \bigl ( f_*T_pM\cap \linsp \{ \xi _p,\widehat {\xi }_p\} \bigr ) =1$ and we can find the vector ${X_2}_p\in T_pM$ such that $f_*{X_2}_p\in \linsp \{ \xi _p,\widehat {\xi }_p\} $ and $\det (f_*{X_1}_p,f_*{X_2}_p,\xi _p)=1$. In this way we define the vector field $X_2$ such that
\begin{equation}
f_*X_2=a_{11}\, \xi +a_{21}\, \widehat {\xi }
\label{fgx2}
\end{equation}
with some functions $a_{11}$ and $a_{21}$, and
\begin{equation}
\det (f_*X_1,f_*X_2,\xi )=1.
\label{unimod}
\end{equation}
Similarly we may define the vector field $\widehat X_2$ such that
\begin{equation}
\widehat f_*\widehat X_2=a_{12}\, \xi +a_{22}\, \widehat {\xi }
\label{fgx2dasz}
\end{equation}
and
\begin{equation}
\det (\widehat f_*\widehat X_1,\widehat f_*\widehat X_2,\widehat {\xi })=1.
\label{unimoddasz}
\end{equation}
From (\ref{unimod}), (\ref{fgx1}) and (\ref{fgx2}) it follows that $a_{21}=-\, \frac{1}{W}$ and
from (\ref{unimoddasz}), (\ref{fgx1dasz}) and (\ref{fgx2dasz}) we obtain $a_{12}=\frac{1}{W}$.
Since, by (\ref{fgx2}), $a_{11}+a_{21}\, \nu (\widehat {\xi })=0$, and by (\ref{fgx2dasz}) $a_{12}\, \widehat {\nu }(\xi )+a_{22}=0$, we have $a_{11}=\frac{A}{W}$ and $a_{22}=-\, \frac{\widehat A}{W}$.
It follows that
\begin{eqnarray}
\widehat f_*\widehat X_1&=&f_*X_1, \nonumber \\
\widehat f_*\widehat X_2&=&\widehat A\, f_*X_2+\frac{1-A\widehat A}{W}\, \xi , \label{dwiebazy} \\
\widehat {\xi }&=&-\,W\, f_*X_2+A\, \xi .\nonumber
\end{eqnarray}
We have
\begin{equation}
\widetilde {\xi }=\lambda \, \xi +f_*Z, \ \ \ \overline {\xi }=\mu \, \widehat {\xi }+\widehat f_*V.
\end{equation}

Let $Z=z^1\, X_1+z^2\, X_2$ and $V=w^1\, \widehat X_1+w^2\, \widehat X_2$. Let $\widetilde W:=\det (\widehat f-f,\widetilde {\xi },\overline {\xi })$.
\begin{eqnarray}
\widetilde W&=&\det (\widehat f-f,\ \lambda \, \xi +z^1\, f_*X_1+z^2\, f_*X_2,\ \mu \, \widehat {\xi }+w^1\, \widehat f_*\widehat X_1+w^2\, \widehat f_*\widehat X_2) \nonumber \\
&=&\det (\widehat f-f,\ \lambda \, \xi +z^2\, f_*X_2,\ \mu \, \widehat {\xi }+w^2\, \widehat f_*\widehat X_2) \nonumber \\
&=&\det \left ( \widehat f-f,\ \lambda \, \xi +z^2\, \left ( \frac{A}{W}\, \xi - \frac{1}{W}\, \widehat {\xi }\right ) ,\ \mu \, \widehat {\xi }+w^2\, \left ( \frac{1}{W}\, \xi -\frac{\widehat A}{W}\, \widehat {\xi }\right ) \right ) \nonumber \\
&=&\det \left ( \widehat f-f,\ \left ( \lambda +z^2\, \frac{A}{W}\right ) \, \xi -\frac{z^2}{W}\, \widehat {\xi },\ \frac{w^2}{W}\, \xi +\left ( \mu  -w^2\, \frac{\widehat A}{W}\right ) \, \widehat {\xi }\right ) \nonumber \\
&=&\left ( \left ( \lambda +z^2\, \frac{A}{W}\right ) \left ( \mu  -w^2\, \frac{\widehat A}{W}\right ) +\frac{z^2\, w^2}{W^2}\right ) \, \det (\widehat f-f,\xi ,\widehat {\xi }) \nonumber \\
&=&\lambda \, \mu \, W + z^2A\mu -w^2\widehat A\lambda +z^2w^2\frac{1-A\widehat A}{W}.\nonumber
\end{eqnarray}

To compute $\widetilde {\nu }(\overline {\xi })$ we have to write $\overline {\xi }$ in the basis $f_*X_1$, $f_*X_2$, $\widetilde {\xi }$.
\begin{eqnarray}
\overline {\xi }&=&\mu \, \widehat {\xi }+w^1\, \widehat f_*\widehat X_1+w^2\, \widehat f_*\widehat X_2 \nonumber \\
&=&\mu \, \left ( -\,W\, f_*X_2+A\, \xi \right ) +w^1\, f_*X_1+w^2\, \left ( \widehat A\, f_*X_2+\frac{1-A\widehat A}{W}\, \xi \right )  \nonumber \\
&=&\left ( \mu A+w^2\, \frac{1-A\widehat A}{W}\right ) \xi +f_*\left ( w^1X_1+\left ( w^2\widehat A-\mu W\right ) X_2\right )
\nonumber \\
&=&\left ( \mu A+w^2\, \frac{1-A\widehat A}{W}\right ) \left ( \frac{1}{\lambda }\, \widetilde {\xi }-\frac{1}{\lambda }\, f_*Z\right ) +f_*\left ( w^1X_1+\left ( w^2\widehat A-\mu W\right ) X_2\right ) .\nonumber
\end{eqnarray}
It follows that
\begin{equation}
\widetilde A:=\widetilde {\nu }(\overline {\xi })=\frac{1}{\lambda }\, \left ( \mu A+w^2\, \frac{1-A\widehat A}{W}\right ) .
\end{equation}
Similarly we obtain
\begin{equation}
\overline A:=\overline {\nu }(\widetilde {\xi })=\frac{1}{\mu }\, \left ( \lambda \widehat A-z^2\, \frac{1-A\widehat A}{W}\right ) .
\end{equation}
Consequently
\begin{equation}
1-\widetilde A\, \overline A=\frac{1-A\widehat A}{\lambda \mu W}\, \widetilde W.
\end{equation}
\medskip

(iv) Since ${\det }_{\widetilde {\theta }}\widetilde h=\frac{1}{\lambda ^4}\, \det _{\theta }h$ and
${\det }_{\overline {\theta }}\overline h=\frac{1}{\mu ^4}\, {\det }_{\widehat {\theta }}\widehat h$ \cite{NS}, we obtain from (iii)
\[ \left ( \frac{1-\widetilde A\, \overline A}{\widetilde W}\right ) ^4\frac{1}{{\det }_{\widetilde {\theta }}\widetilde h\,
{\det }_{\overline {\theta }}\overline h}=\left ( \frac{1-A\widehat A}{\lambda \mu W}\right ) ^4\frac{\lambda ^4\mu ^4}{\det _{\theta }h\, {\det }_{\widehat {\theta }}\widehat h}=\left ( \frac{1-A\widehat A}{W}\right ) ^4\frac{1}{\det _{\theta }h\, {\det }_{\widehat {\theta }}\widehat h}. \]
\bigskip

From (iv) of Proposition \ref{funkcja} it follows that
\begin{equation}
\psi (f,\widehat f):=\left ( \frac{1-\nu (\widehat {\xi })\,
\widehat {\nu }(\xi )}{\det (\widehat f-f,\xi ,\widehat {\xi })}\right ) ^4\frac{1}{{\det }_{\theta }h\ {\det }_{\widehat {\theta }}\widehat h}
\end{equation}
is a well defined function on $M$.
\bigskip

Throughout the paper we will make some assumption about the rank of the spherical representation of $\widehat f-f$. The following lemma explains the technical significance of this assumption: the forms $\omega ^2_{\ 1}$, $\omega ^3_{\ 1}$ constitute a local frame of $T^*M$.
\begin{lemma}
\emph{(cf \cite{Sh} page 6 in the metric case)} Let $\varphi :M\to GL(3,\R )$. For $p\in M$ we denote by ${v_1}_p$, ${v_2}_p$, ${v_3}_p$ the columns of the matrix $\varphi (p)$. We consider the mappings $v_1:M\to \R ^3\setminus \{ \mathbf{0} \} $ and $\pi \circ v_1:M\to \Proj ^2(\R )$, where $\pi :\R ^3\setminus \{ \mathbf{0} \} \to \Proj ^2(\R )$ denotes the canonical projection. The forms $\omega ^i_{\ 1}$ are defined by the equality
\begin{equation}
dv_1=\omega ^1_{\ 1}\,v_1+\omega ^2_{\ 1}\, v_2+\omega ^3_{\ 1}v_3.
\label{dv1}
\end{equation}
At each point of $M$ the following conditions are equivalent:

\emph{(i)} $\mbox{\emph{rank}}(\pi \circ v_1)=2$,

\emph{(ii)} $\omega ^2_{\ 1}\wedge \omega ^3_{\ 1}\ne 0$.
\label{rankspherrep2}
\end{lemma}
\emph{Proof.} Let $(v_i^1,v_i^2,v_i^3)$ be the coordinates of $v_i$. Assume for example that $v_1^3\ne 0$. Then on $\Proj ^2(\R )$ we use the chart $(t^1:t^2:t^3)\mapsto \left ( \frac{t^1}{t^3},\frac{t^2}{t^3}\right ) $. The composition of $\pi \circ v_1$ with this chart is $\left ( \frac{v_1^1}{v_1^3},\frac{v_1^2}{v_1^3}\right ) $ and its rank equals two if and only if $d\left ( \frac{v_1^1}{v_1^3}\right ) \wedge d\left ( \frac{v_1^2}{v_1^3}\right ) \ne 0$.
Let $Z(p)$ be the inverse matrix of $\varphi (p)$ and let $Z=(z_{ij})$.
Using (\ref{dv1}) we easily obtain $d\left ( \frac{v_1^1}{v_1^3}\right ) =\frac{\det \varphi }{(v_1^3)^2}\, \left ( z_{32}\omega ^2_{\ 1}-z_{22}\omega ^3_{\ 1}\right ) $, $d\left ( \frac{v_1^2}{v_1^3}\right ) =\frac{\det \varphi }{(v_1^3)^2}\, \left ( -\, z_{31}\omega ^2_{\ 1}+z_{21}\omega ^3_{\ 1}\right ) $ and
\begin{eqnarray}
d\left ( \frac{v_1^1}{v_1^3}\right ) \wedge d\left ( \frac{v_1^2}{v_1^3}\right ) &=&
\frac{(\det \varphi )^2}{(v_1^3)^4}\, \left | \begin{array}{cc} z_{21} & z_{22} \\ z_{31} & z_{32} \end{array} \right | \, \omega ^2_{\ 1}\wedge \omega ^3_{\ 1} \nonumber \\
&=&\frac{(\det \varphi )^2}{(v_1^3)^4}\, \det Z\, v_1^3\, \omega ^2_{\ 1}\wedge \omega ^3_{\ 1}=\frac{\det \varphi }{(v_1^3)^3}\, \omega ^2_{\ 1}\wedge \omega ^3_{\ 1}, \nonumber
\end{eqnarray}
hence $d\left ( \frac{v_1^1}{v_1^3}\right ) \wedge d\left ( \frac{v_1^2}{v_1^3}\right ) \ne 0$ is equivalent to $\omega ^2_{\ 1}\wedge \omega ^3_{\ 1}\ne 0$.
If at some point $v_1^3=0$, then we have to use another chart and one of the equalities
$d\left ( \frac{v_1^1}{v_1^2}\right ) \wedge d\left ( \frac{v_1^3}{v_1^2}\right ) =-\,
\frac{\det \varphi }{(v_1^2)^3}\, \omega ^2_{\ 1}\wedge \omega ^3_{\ 1}$,
$d\left ( \frac{v_1^2}{v_1^1}\right ) \wedge d\left ( \frac{v_1^3}{v_1^1}\right ) =
\frac{\det \varphi }{(v_1^1)^3}\, \omega ^2_{\ 1}\wedge \omega ^3_{\ 1}$.

\begin{theorem}
Let $f$ and $\widehat f$ be as in Proposition \emph{\ref{funkcja}}. Assume that the spherical representation of $\widehat f-f$,
$M\ni p\mapsto \pi (\widehat f(p)-f(p))\in \Proj ^2(\R )$,
has rank $2$ at every point of $M$ .
\smallskip

Then:

\emph{(i)} $f_*T_pM\ne \widehat f_*T_pM$ for every $p\in M$,

\emph{(ii)} $\widehat f(p)-f(p)$ is not an asymptotic vector,

\emph{(iii)} the affine fundamental forms $h$ and $\widehat h$ are
conformal to each other if and only if $\psi (f,\widehat f)=1$.
\label{wkwhconfdaszh}
\end{theorem}

\emph{Proof.} We choose transversal fields $\xi $ and $\widehat {\xi }$ satisfying $\det (\widehat f-f, \xi ,\widehat {\xi })\ne 0$. We retain the notation of Proposition \ref{funkcja} and Lemma \ref{rankspherrep2}.
We take
\[ v_1=\widehat v_1=\widehat f-f, \ \ \ v_2=f_*X_2, \ \ \ \widehat v_2=\widehat f_*\widehat X_2, \ \ \
v_3=\xi \ \ \mbox{and} \ \ \widehat v_3=\widehat {\xi }.\]
Together with $f$ and $\widehat f$ we consider moving frames $F$ and $\widehat F$ from $M$ to $ASL(3,\R )$,
\[ F=\left ( \begin{array}{cc} 1 & 0 \\ f & (v_1,v_2,v_3) \end{array} \right ) , \ \ \
\widehat F= \left ( \begin{array}{cc} 1 & 0 \\ \widehat f & (\widehat v_1,\widehat v_2,\widehat v_3) \end{array} \right ) ,\]
We can now rewrite (\ref{fgx1}) and (\ref{dwiebazy}) as $\widehat F=Fa$ with
\[ a=\left ( \begin{array}{cccc} 1 & 0 & 0 & 0 \\ 1 & 1 & 0 & 0 \\ 0 & 0 & \widehat A & -\, W \\ 0 & 0 & \frac{1-A\widehat A}{W} & A \end{array} \right ) . \]
The pull-back of the Maurer-Cartan form $\vartheta $ on $ASL(3,\R )$ by $F$ is
\[ F^*\vartheta =F^{-1}\, dF=\left ( \begin{array}{cccc} 0 & 0 & 0 & 0 \\ \vartheta ^1 & \omega ^1_{\ 1} & \omega ^1_{\ 2} & \omega ^1_{\ 3} \\ \vartheta ^2 & \omega ^2_{\ 1} & \omega ^2_{\ 2} & \omega ^2_{\ 3}
\\ \vartheta ^3 & \omega ^3_{\ 1} & \omega ^3_{\ 2} & \omega ^3_{\ 3} \end{array} \right ) .\]
Then
\begin{eqnarray}
df&=&\vartheta ^1\, v_1+\vartheta ^2\, v_2+\vartheta ^3\, v_3, \nonumber \\
dv_1&=&\omega ^1_{\ 1}\, v_1+\omega ^2_{\ 1}\, v_2+\omega ^3_{\ 1}\, v_3, \nonumber \\
dv_2&=&\omega ^1_{\ 2}\, v_1+\omega ^2_{\ 2}\, v_2+\omega ^3_{\ 2}\, v_3, \nonumber \\
dv_3&=&\omega ^1_{\ 3}\, v_1+\omega ^2_{\ 3}\, v_2+\omega ^3_{\ 3}\, v_3. \nonumber
\end{eqnarray}
Since $d\circ d =0$, the $1$-forms $\vartheta ^i$ and $\omega ^j_{\ k}$ satisfy the structure equations
\begin{equation}
d\vartheta ^s=-\, \sum _{k=1}^3\omega ^s_{\ k}\wedge \vartheta ^k, \ \ s=1,2,3;
\label{streqtheta}
\end{equation}
and
\begin{equation}
d\omega ^i_{\ j}=-\, \sum _{k=1}^3\omega ^i_{\ k}\wedge \omega ^k_{\ j}, \ \ \ i,j=1,2,3.
\end{equation}
Similar equalities one can write for the dashed $1$-forms $\widehat {\vartheta }^i$ and $\widehat {\omega }^j_{\ k}$.
\bigskip

From
\begin{equation}
{\widehat F}^{-1}d\widehat F=a^{-1}\, (F^{-1}dF)\, a+a^{-1}da
\label{zmformkoneks}
\end{equation}
we obtain
\begin{eqnarray}
\widehat {\vartheta }^2&=&A\, \vartheta ^2+W\, \vartheta ^3+A\, \omega ^2_{\ 1}+W\, \omega ^3_{\ 1}, \label{dasztheta2} \\
\widehat {\vartheta }^3&=&-\, \frac{1-A\widehat A}{W}\, \vartheta ^2+\widehat A\, \vartheta ^3-\frac{1-A\widehat A}{W}\, \omega ^2_{\ 1}+\widehat A\, \omega ^3_{\ 1} \label{dasztheta3}
\end{eqnarray}
Since the frames $(v_1,v_2,v_3)$ and $(\widehat v_1,\widehat v_2,\widehat v_3)$ are adapted to $f$ and $\widehat f$ respectively, we have $\vartheta ^3=0$ and
$\widehat {\vartheta }^3=0$. From (\ref{dasztheta3}) we obtain
\begin{equation}
0=-\, \frac{1-A\widehat A}{W}\left ( \vartheta ^2+\omega ^2_{\ 1}\right ) +\widehat A\, \omega ^3_{\ 1}.
\label{dasztheta3-2}
\end{equation}
\bigskip
Suppose that $1-A\widehat A=0$. Then (\ref{dasztheta3-2}) and $\omega ^2_{\ 1}\wedge \omega ^3_{\ 1}\ne 0$ imply $\widehat A=0$, which contradicts $1-A\widehat A=0$. Therefore $1-A\widehat A\ne 0$ and from
(\ref{dwiebazy}) we obtain (i).
\smallskip

From (\ref{dasztheta3-2}) and (\ref{dasztheta2}) it follows that
\begin{eqnarray}
\vartheta ^2&=&-\, \omega ^2_{\ 1}+\frac{\widehat AW}{1-A\widehat A}\, \omega ^3_{\ 1}, \label{theta2} \\
\widehat {\vartheta }^2&=&\frac{W}{1-A\widehat A}\, \omega ^3_{\ 1}. \label{dasztheta2-2}
\end{eqnarray}

From (\ref{zmformkoneks}) we obtain also
\begin{eqnarray}
\widehat {\omega }^2_{\ 1}&=&A\, \omega ^2_{\ 1}+W\, \omega ^3_{\ 1}, \label{daszomega21} \\
\widehat {\omega }^3_{\ 1}&=&-\, \frac{1-A\widehat A}{W}\, \omega ^2_{\ 1}+\widehat A\, \omega ^3_{\ 1}. \label{daszomega31}
\end{eqnarray}
Comparing (\ref{theta2}) with (\ref{daszomega31}) yields
\begin{equation}
\vartheta ^2=\frac{W}{1-A\widehat A}\, \widehat {\omega }^3_{\ 1}. \label{theta2-3}
\end{equation}

Our next goal is to check that $X_1$ and $\widehat X_1$ are at each point linearly independent. We only need to show that $\vartheta ^2\wedge \widehat {\vartheta }^2\ne 0$ and it suffices to use (\ref{theta2}) and (\ref{dasztheta2-2}) to obtain
\[ \vartheta ^2\wedge \widehat {\vartheta }^2=\frac{-\, W}{1-A\widehat A}\, \omega ^2_{\ 1}\wedge \omega ^3_{\ 1}.\]

We may now find the matrices of $h$ and $\widehat h$ in the basis $X_1$, $\widehat X_1$. Since for $k\in \{ 1,2\} $ $h(Y,X_k)=\omega ^3_{\ k}(Y)$ and $\widehat h(Y,\widehat X_k)=\widehat {\omega }^3_{\ k}(Y)$, we obtain from (\ref{dasztheta2-2}) and (\ref{theta2-3})
\begin{equation}
h(\widehat X_1,X_1)=0 \ \ \ \mbox{and} \ \ \ \widehat h(X_1,\widehat X_1)=0.
\end{equation}
It follows that $h(X_1,X_1)\ne 0$ and $\widehat h(\widehat X_1,\widehat X_1)\ne 0$, for otherwise $f$ or $\widehat f$ would be degenerate. We thus get (ii).
\smallskip

Let $h_{ij}=h(X_i,X_j)$ and $\widehat h_{ij}=\widehat h(\widehat X_i,\widehat X_j)$. Let $\widehat X_1=c_{11}\, X_1+c_{21}\, X_2$. Here
\[ c_{21}=\vartheta ^2(\widehat X_1)=\frac{W}{1-A\widehat A}\, \widehat {\omega }^3_{\ 1}(\widehat X_1)=\frac{W}{1-A\widehat A}\, \widehat h_{11}\]
and consequently
\[ h(\widehat X_1,\widehat X_1)=h(\widehat X_1,c_{11}\, X_1+c_{21}\, X_2)=c_{21}\, h(\widehat X_1,X_2)=\frac{W\, \widehat h_{11}}{1-A\widehat A}\, \omega ^3_{\ 2}(\widehat X_1).\]
In a similar way we obtain
\[ \widehat h(X_1,X_1)=\frac{W\, h_{11}}{1-A\widehat A}\, \widehat {\omega }^3_{\ 2}(X_1).\]
Since $h(X_1,\widehat X_1)=0$ and $\widehat h(X_1,\widehat
X_1)=0$, the affine fundamental form $\widehat h$ is conformal to
$h$ if and only if there exists a function $\lambda $ such that
$\widehat h(X_1,X_1)=\lambda \, h(X_1,X_1)$ and $\widehat
h(\widehat X_1,\widehat X_1)=\lambda \, h(\widehat X_1,\widehat
X_1)$, which is equivalent to
\begin{equation}
\left | \begin{array}{cc} \frac{W\, h_{11}}{1-A\widehat A}\, \widehat {\omega }^3_{\ 2}(X_1) & h_{11} \\
\widehat h_{11} & \frac{W\, \widehat h_{11}}{1-A\widehat A}\,
\omega ^3_{\ 2}(\widehat X_1) \end{array} \right | =0.
\label{wkw1}
\end{equation}
The left-hand side of (\ref{wkw1}) equals $0$ if and only if
\begin{equation}
\left ( \frac{W}{1-A\widehat A}\right ) ^2\, \widehat {\omega }
^3_{\ 2}(X_1)\, \omega ^3_{\ 2}(\widehat X_1)=1, \label{wkw2}
\end{equation}
because $h_{11}\, \widehat h_{11}\ne 0$. Let $H:=\det _{\theta }h$
and $\widehat H:=\det _{\widehat {\theta }}\widehat h$. We have
\begin{eqnarray}
& &H\, \vartheta ^1\wedge \vartheta ^2(X_1,X_2)=H=h_{11}\,
h_{22}-h_{12}\, h_{12} \nonumber \\
&=&\omega ^3_{\ 1}(X_1)\, \omega ^3_{\ 2}(X_2)-\omega ^3_{\
1}(X_2)\, \omega ^3_{\ 2}(X_1)=\omega ^3_{\ 1}\wedge \omega ^3_{\
2}(X_1,X_2),
\end{eqnarray}
hence
\begin{equation}
\omega ^3_{\ 1}\wedge \omega ^3_{\ 2}=H\, \vartheta ^1\wedge
\vartheta ^2. \label{omega31wedgeomega32}
\end{equation}
Similarly
\begin{equation}
\widehat {\omega }^3_{\ 1}\wedge \widehat {\omega }^3_{\
2}=\widehat H\, \widehat {\vartheta }^1\wedge \widehat {\vartheta
}^2. \label{daszomega31wedgedaszomega32}
\end{equation}
Using (\ref{dasztheta2-2}) and (\ref{theta2-3}) we obtain
\begin{equation}
\omega ^3_{\ 1}\wedge \omega ^3_{\ 2}(\widehat X_1,\widehat
X_2)=\frac{1-A\widehat A}{W}\, \widehat {\vartheta }^2\wedge
\omega ^3_{\ 2}(\widehat X_1,\widehat X_2)=-\, \frac{1-A\widehat
A}{W}\, \omega ^3_{\ 2}(\widehat X_1) \label{omega32nadaszX1}
\end{equation}
\begin{equation}
\widehat {\omega }^3_{\ 1}\wedge \widehat {\omega }^3_{\
2}(X_1,X_2)=\frac{1-A\widehat A}{W}\, \vartheta ^2\wedge \widehat
{\omega }^3_{\ 2}(X_1,X_2)=-\, \frac{1-A\widehat A}{W}\, \widehat
{\omega }^3_{\ 2}(X_1) \label{daszomega32naX1}
\end{equation}
Combining (\ref{omega32nadaszX1}) with (\ref{omega31wedgeomega32})
and (\ref{daszomega32naX1}) with
(\ref{daszomega31wedgedaszomega32}) gives
\begin{equation}
\omega ^3_{\ 2}(\widehat X_1)=\frac{-\, WH}{1-A\widehat A}\,
\vartheta ^1\wedge \vartheta ^2(\widehat X_1,\widehat X_2)
\end{equation}
and
\begin{equation}
\widehat {\omega }^3_{\ 2}(X_1)=\frac{-\, W\widehat H}{1-A\widehat
A}\, \widehat {\vartheta }^1\wedge \widehat {\vartheta
}^2(X_1,X_2).
\end{equation}
Condition (\ref{wkw2}) now becomes
\begin{equation}
\left ( \frac{W}{1-A\widehat A}\right ) ^4\, H\widehat H\,
\vartheta ^1\wedge \vartheta ^2(\widehat X_1,\widehat X_2)\,
\widehat {\vartheta }^1\wedge \widehat {\vartheta }^2(X_1,X_2)=1. \label{wkw3}
\end{equation}
But $\vartheta ^1\wedge \vartheta ^2(\widehat X_1,\widehat X_2)\,
\widehat {\vartheta }^1\wedge \widehat {\vartheta }^2(X_1,X_2)=1$,
because the matrix $(\widehat {\vartheta }^i(X_j))$ is the inverse
of $(\vartheta ^k(\widehat X_l))$. We thus get (iii).
\bigskip

As a supplement we give here another similar criterion, applicable when we want to use parallel transversal fields $\xi $ and $\widehat {\xi }$. The equality in (iii) corresponds to (3.22) in \cite{ChT}.
\medskip

\begin{theorem}
\label{hconfdaszh-parallelxi}
Let $f$, $\widehat f$ be as in Proposition \emph{\ref{funkcja}} and let $X_1$, $\widehat X_1$ satisfy \emph{(\ref{fgx1})} and \emph{(\ref{fgx1dasz})}. Assume that $\xi $ and $\widehat {\xi }$, transversal fields for $f$ and $\widehat f$ respectively, are parallel.
\smallskip

We choose arbitrary $X_2$ such that $X_1$, $X_2$ is a local frame unimodular with respect to $\theta _{\xi }$. Let $\widehat X_2$ be defined by the following two conditions: for every $p\in M$ $\widehat f_*(T_pM)\cap \linsp \{ f_*(X_{2\, p}), \ \xi _p\} =\R \, \widehat f_*(\widehat X_{2\, p})$ and $\widehat {\theta }_{\widehat {\xi }}(\widehat X_1,\widehat X_2)=1$. 
\smallskip

Then:
\smallskip

\emph{(i)} $\widehat f_*(\widehat X_2)=\lambda \, f_*(X_2)+\beta \, \xi , \ \ \ \ \widehat {\xi }=\frac{1}{\lambda }\, \xi $ for some functions $\lambda $, $\beta $;
\smallskip

\emph{(ii)} $\lambda $, $\beta $ do not depend on $X_2$ ($\widehat X_2$ does);
\smallskip

\emph{(iii)} if the spherical representation 
$\pi \circ (\widehat f-f)$ of $\widehat f-f$
has rank $2$ at every point of $M$, then affine fundamental forms $h$ and $\widehat h$ are proportional if and only if 
$\det _{\theta }h\cdot \det _{\widehat {\theta }}\widehat h=\beta ^{\, 4}$.
\end{theorem} 
\emph{Proof.} By assumption, $\widehat f_*(\widehat X_2)=\lambda \, f_*(X_2)+\beta \, \xi $ and $\widehat {\xi }=\mu \, \xi $ for some functions $\lambda $, $\mu $ and $\beta $. From $\widehat {\theta }_{\widehat {\xi }}(\widehat X_1,\widehat X_2)=1$ we obtain $\mu \cdot \lambda =1$ and (i) follows. 

If we replace $X_2$ by $X_2+t\, X_1$, then $\widehat X_2$ should be replaced by $\widehat X_2+\lambda t\, \widehat X_1$. We have then $\widehat f_*(\widehat X_2+\lambda t\, \widehat X_1)=\lambda f_*(X_2+t\, X_1)+\beta \, \xi $.

Note that $\beta \ne 0$, because $\beta =0$ would imply $\omega ^3_{\ 1}=0$, which contradicts the non-degeneracy of $f$.
\smallskip

Proof of (iii) is similar to the proof of (iii) in Theorem \ref{wkwhconfdaszh}. 
We have now $\widehat F=F\, a$ with $a=\left ( \begin{array}{cccc} 1 & 0 & 0 & 0 \\ 1 & 1 & 0 & 0 \\ 0 & 0 & \lambda & 0 \\ 0 & 0 & \beta & \frac{1}{\lambda }\end{array} \right ) $ and from (\ref{zmformkoneks}) we obtain in particular 
\begin{align}
\omega ^3_{\ 1}&=\beta \, \widehat {\vartheta }^2, \\
\widehat {\omega }^3_{\ 1}&=\beta \, \vartheta ^2, \\
\vartheta ^2\wedge \widehat {\vartheta }^2&=\frac{-\, 1}{\beta }\, \omega ^2_{\ 1}\wedge \omega ^3_{\ 1}\ne 0.
\end{align}
The rest of the proof runs as before, with $\frac{1-A\widehat A}{W}$ replaced by $\beta $.

We may also compute $\psi (f,\widehat f)$ using $\xi $ and $\widehat {\widehat {\xi }}=f_*(X_2)=\frac{1}{\lambda }\, \widehat f_*(\widehat X_2)-\beta \, \widehat {\xi }$ as linearly independent transversal fields for $f$and $\widehat f$ respectively, and apply Theorem \ref{wkwhconfdaszh}. Note that $\nu (\widehat {\widehat {\xi }})=0$ and $\det _{\widehat {\widehat {\theta }}}\widehat {\widehat h}=\frac{1}{\beta ^4}\, \det _{\widehat {\theta }}\widehat h$.

\section{Chern--Terng theorem}

\begin{theorem} \emph{\cite{ChT}}
Let $\dim M=2$ and $f, \widehat f:M\to \R ^3$, be a pair of
non-degenerate immersions, satisfying the following conditions:
\smallskip

\emph{(i)} for every $p\in M$: $f(p)\ne \widehat f(p)$,
the vector $\widehat f(p)-f(p)$ is tangent to $f(M)$ at $f(p)$ and is tangent to $\widehat f(M)$ at $\widehat f(p)$,
\smallskip

\emph{(ii)} the affine fundamental forms of $f$ and $\widehat f$ are conformal to each other,

\emph{(iii)} the affine normals of both surfaces at corresponding points $f(p)$ and $\widehat f(p)$ are parallel.
\smallskip

Then the surfaces are both affine minimal.
\end{theorem}
\emph{Proof.} We give here a proof which in some details will be different from that in \cite{ChT}, because we want to use local frames with the last vector field equal to corresponding affine normal vector field. 
\smallskip

At first we consider the set of points where the rank of the spherical representation of $\widehat f-f$ equals 2. 
We use the same local frame as in Theorem \ref{hconfdaszh-parallelxi}.
From assumption (ii) and from Theorem \ref{hconfdaszh-parallelxi} we have $H\cdot \widehat H=\beta ^4$. Since $\xi $ and $\widehat {\xi }$ are affine normal vector fields, $|H|=1$ and $|\widehat H|=1$. It follows that $|\beta |=1$. If we replace $\widehat {\xi }$ by $-\, \widehat {\xi }$, then $\widehat X_2$ should be replaced by $-\, \widehat X_2$, $\lambda $ by $-\, \lambda $ and $\beta $ by $-\, \beta $. Therefore without loss of generality we may assume that $\beta =1$. Moreover, $H=\widehat H=:\varepsilon _h$, because $H\cdot \widehat H>0$.
\smallskip

From (\ref{zmformkoneks}) we obtain $\widehat {\vartheta }^3=-\, \beta \vartheta ^2-\beta \, \omega ^2_{\ 1}+\lambda \, \omega ^3_{\ 1}$ and $\widehat {\omega }^3_{\ 3}=-\, \frac{\beta }{\lambda }\, \omega ^2_{\ 3}-\frac{d\lambda }{\lambda }$. Then $\widehat {\vartheta }^3=0$, $\widehat {\omega }^3_{\ 3}=0$ together with $\beta =1$ give 
\begin{align}
\vartheta ^2+\omega ^2_{\ 1}&=\lambda \, \omega ^3_{\ 1}, \\
d\lambda +\omega ^2_{\ 3}&=0,
\end{align}
which corresponds to $\gamma =0$ and $\beta =0$ in (3.8) of \cite{ChT}. We will next assume that $\varepsilon _h+\lambda ^2\ne 0$ and prove the equality corresponding to $\alpha =0$, that is 
\begin{equation}
\vartheta ^1+\omega ^1_{\ 1}=-\, \lambda \, \omega ^3_{\ 2}.
\end{equation} 
Application of (\ref{zmformkoneks}) gives
\begin{align*}
\widehat {\vartheta }^1&=\vartheta ^1+\omega ^1_{\ 1}, \\
\widehat {\omega }^3_{\ 2}&=-\, \lambda \beta \, \omega ^2_{\ 2}+\lambda ^2\, \omega ^3_{\ 2}-\, \beta ( \beta \, \omega ^2_{\ 3}+d\lambda )+\lambda \, d\beta =\lambda \, \omega ^1_{\ 1}+\lambda ^2\, \omega ^3_{\ 2}.
\end{align*}
Let $\varphi =\vartheta ^1+\omega ^1_{\ 1}+\lambda \, \omega ^3_{\ 2}$.
We have
\[ 
\widehat {\vartheta }^1=\varphi -\lambda \, \omega ^3_{\ 2}, \ \ \ \
\widehat {\vartheta }^2=\omega ^3_{\ 1}, \ \ \ \
\widehat {\omega }^3_{\ 1}=\vartheta ^2, \ \ \ \
\widehat {\omega }^3_{\ 2}=\lambda \, (\varphi - \vartheta ^1).
\]
Then
\begin{align*}
0&=\widehat {\omega }^3_{\ 1}\wedge \widehat {\vartheta }^1+\widehat {\omega }^3_{\ 2}\wedge \widehat {\vartheta }^2
=\vartheta ^2\wedge (\varphi -\lambda \, \omega ^3_{\ 2})+\lambda \, (\varphi - \vartheta ^1)\wedge \omega ^3_{\ 1} \\
&=(\vartheta ^2-\lambda \, \omega ^3_{\ 1})\wedge \varphi +\lambda (\omega ^3_{\ 1}\wedge \vartheta ^1+\omega ^3_{\ 2}\wedge \vartheta ^2)=(\vartheta ^2-\lambda \, \omega ^3_{\ 1})\wedge \varphi 
\end{align*}
and
\begin{align*}
0&=\widehat {\omega }^3_{\ 1}\wedge \widehat {\omega }^3_{\ 2}-\varepsilon _h\, \widehat {\vartheta }^1\wedge \widehat {\vartheta }^2=\vartheta ^2\wedge \lambda \, (\varphi - \vartheta ^1)-\varepsilon _h\, (\varphi -\lambda \, \omega ^3_{\ 2})\wedge \omega ^3_{\ 1} \\
&=(\lambda \, \vartheta ^2+\varepsilon _h\, \omega ^3_{ 1})\wedge \varphi +\varepsilon _h\, \lambda (\varepsilon _h\, \vartheta ^1\wedge \vartheta ^2-\omega ^3_{\ 1}\wedge \omega ^3_{\ 2})=(\lambda \, \vartheta ^2+\varepsilon _h\, \omega ^3_{ 1})\wedge \varphi .
\end{align*}
If $\varepsilon _h+\lambda ^2\ne 0$, then the 1-forms $\vartheta ^2-\lambda \, \omega ^3_{\ 1}$ and $\lambda \, \vartheta ^2+\varepsilon _h\, \omega ^3_{ 1}$ are linearly independent, because $(\vartheta ^2-\lambda \, \omega ^3_{\ 1})\wedge (\lambda \, \vartheta ^2+\varepsilon _h\, \omega ^3_{ 1})=(\varepsilon _h+\lambda ^2)\, \vartheta ^2\wedge \omega ^3_{\ 1}\ne 0$ (recall that in the considered case $h_{11}\ne 0$). Consequently the equalities 
$(\vartheta ^2-\lambda \, \omega ^3_{\ 1})\wedge \varphi =0$ and $(\lambda \, \vartheta ^2+\varepsilon _h\, \omega ^3_{ 1})\wedge \varphi =0$ imply $\varphi =0$.
\smallskip

It follows that
\begin{align*}
0&=d\widehat {\omega }^3_{\ 3}=-\, \widehat {\omega }^3_{\ 1}\wedge \widehat {\omega }^1_{\ 3}-\widehat {\omega }^3_{\ 2}\wedge \widehat {\omega }^2_{\ 3}=-\, \vartheta ^2\wedge \frac{1}{\lambda }\, \omega ^1_{\ 3}+\lambda \, \vartheta ^1\wedge \frac{1}{\lambda ^2}\, \omega ^2_{\ 3} \\
&=\frac{1}{\lambda }\, (-\, \vartheta ^2\wedge \omega ^1_{\ 3}+\vartheta ^1\wedge \omega ^2_{\ 3}),
\end{align*}
which implies $\mbox{tr}\, S=0$, and 
\begin{align*}
0&=d\omega ^3_{\ 3}=-\, \omega ^3_{\ 1}\wedge \omega ^1_{\ 3}-\omega ^3_{\ 2}\wedge \omega ^2_{\ 3}=-\, \widehat {\vartheta }^2\wedge \lambda \, \widehat {\omega }^1_{\ 3}+\frac{1}{\lambda }\, \widehat {\vartheta }^1\wedge \lambda ^2\, \widehat {\omega }^2_{\ 3} \\
&=\lambda \, (-\, \widehat {\vartheta }^2\wedge \widehat {\omega }^1_{\ 3}+\widehat {\vartheta }^1\wedge \widehat {\omega }^2_{\ 3}),
\end{align*}
hence also $\mbox{tr}\, \widehat S=0$.
\smallskip

We thus get $\mbox{tr}\, S=0$ and $\mbox{tr}\, \widehat S=0$ on the set of points where $\mbox{rank}\,( \pi \circ (\widehat f-f))=2$ and $\varepsilon _h+\lambda ^2\ne 0$, and also on its closure, by continuity.
\medskip

Assume now that $\varepsilon _h+\lambda ^2=0$ on some open set, contained in the set where $\mbox{rank}\,( \pi \circ (\widehat f-f))=2$ holds. In this case $d\lambda =0$, hence $\omega ^2_{\ 3}=0$. We have
\begin{align*}
0&=d\omega ^2_{\ 3}=-\, \omega ^2_{\ 1}\wedge \omega ^1_{\ 3}-\omega ^2_{\ 2}\wedge \omega ^2_{\ 3}=-\, \omega ^2_{\ 1}\wedge \omega ^1_{\ 3}, \\
0&=d\omega ^3_{\ 3}=-\, \omega ^3_{\ 1}\wedge \omega ^1_{\ 3}-\omega ^3_{\ 2}\wedge \omega ^2_{\ 3}=-\, \omega ^3_{\ 1}\wedge \omega ^1_{\ 3},
\end{align*}
and it follows that 
\[ \vartheta ^2\wedge \omega ^1_{\ 3}-\vartheta ^1\wedge \omega ^2_{\ 3}=\vartheta ^2\wedge \omega ^1_{\ 3}=
(-\, \omega ^2_{\ 1}+\lambda \, \omega ^3_{\ 1})\wedge \omega ^1_{\ 3}=0\]
and 
\[ \widehat {\vartheta }^2\wedge \widehat {\omega }^1_{\ 3}-\widehat {\vartheta }^1\wedge \widehat {\omega }^2_{\ 3}=\omega ^3_{\ 1}\wedge \frac{1}{\lambda }\, \omega ^1_{\ 3}-\widehat {\vartheta }^1\wedge \frac{1}{\lambda ^2}\, \omega ^2_{\ 3}=0.\]
\medskip

Finally, we consider the interior of the set where 
$\mbox{rank}\, (\pi \circ (\widehat f-f))<2$. Since $\omega ^3_{\ 1}\ne 0$, $\mbox{rank}\, (\pi \circ (\widehat f-f))\ne 0$. By Lemma \ref{rankspherrep2}, $\omega ^2_{\ 1}\wedge \omega ^3_{\ 1}=0$. We will show that also in this case proportionality of $h$ and $\widehat h$ implies $|\beta |=1$, $d\lambda =0$ and $\omega ^2_{\ 3}=0$ as in the preceding case.
\smallskip

From (\ref{zmformkoneks}) we get $\beta \, \vartheta ^2+\beta \, \omega ^2_{\ 1}=\lambda \, \omega ^3_{\ 1}$. Then 
$\omega ^2_{\ 1}\wedge \omega ^3_{\ 1}=0$ and $\beta \ne 0$ imply $\omega ^3_{\ 1}\wedge \vartheta ^2=0$, in particular 
$h_{11}=h(X_1,X_1)=\omega ^3_{\ 1}\wedge \vartheta ^2(X_1,X_2)=0$. Since $\xi $ is an affine normal vector field, 
\[
1=|H|=|h_{11}\, h_{22}-h_{12}^2|=|h_{12}|^2,
\]
hence $h(X_1,X_2)=h_{12}=\varepsilon _1\in \{ 1,-1\} $ and we see that 
\begin{equation}
\omega ^3_{\ 1}=\varepsilon _1\, \vartheta ^2 \ \ \ \mbox{and} \ \ \ \omega ^3_{\ 2}=\varepsilon _1\, \vartheta ^1+h_{22}\, \vartheta ^2.
\label{omega31omega32CT}
\end{equation}
From (\ref{omega31omega32CT}) and (\ref{zmformkoneks}) we have $\varepsilon _1\, \vartheta ^2=\omega ^3_{\ 1}=\beta \, \widehat {\vartheta }^2$ and it follows that $\widehat X_1=c_{11}\, X_1$ for some function $c_{11}$. Then
\[
\widehat h_{11}=\widehat h(\widehat X_1,\widehat X_1)=c_{11}^2\, \widehat h(X_1,X_1)=0,
\]
because $h_{11}=0$ and $\widehat h$ is proportional to $h$. Now from $|\widehat H|=1$ we easily obtain 
\begin{equation}
\widehat {\omega }^3_{\ 1}=\varepsilon _2\, \widehat {\vartheta }^2
\end{equation}
and consequently
\[ \vartheta ^2=\varepsilon _1\, \omega ^3_{\ 1}=\varepsilon _1\, \beta \, \widehat {\vartheta }^2=\varepsilon _1\, \beta \, \varepsilon _2\, \widehat {\omega }^3_{\ 1}=\varepsilon _1\, \varepsilon _2\, \beta ^2\, \vartheta ^2, \]
hence
$|\beta |=1$ and $\varepsilon _1=\varepsilon _2$. Without loss of generality we may assume that $\beta =1$. 
\smallskip

Differentiating both sides of $\omega ^3_{\ 1}=\varepsilon _1\, \vartheta ^2$, using fundamental equations and the equality $\omega ^2_{\ 2}=-\, \omega ^1_{\ 1}$ we obtain
\begin{align}
\omega ^3_{\ 1}\wedge \omega ^1_{\ 1}+\omega ^3_{\ 2}\wedge \omega ^2_{\ 1}&=\varepsilon _1\, \omega ^2_{\ 1}\wedge \vartheta ^1+\varepsilon _1\, \omega ^2_{\ 2}\wedge \vartheta ^2, \nonumber \\
\omega ^3_{\ 2}\wedge \omega ^2_{\ 1}&=\varepsilon _1\, \omega ^2_{\ 1}\wedge \vartheta ^1. \label{o32wo21-1}
\end{align}
We have also
\begin{align}
\omega ^2_{\ 1}\wedge \vartheta ^2&=\omega ^2_{\ 1}\wedge ( -\, \omega ^2_{\ 1}+\lambda \, \omega ^3_{\ 1})=0, 
\label{o21wth2} \\
\omega ^3_{\ 2}\wedge \omega ^2_{\ 1}&=(\varepsilon _1\, \vartheta ^1+h_{22}\, \vartheta ^2)\wedge \omega ^2_{\ 1}=-\, \varepsilon _1\, \omega ^2_{\ 1}\wedge \vartheta ^1. \label{o32wo21-2}
\end{align}
Comparing (\ref{o32wo21-1}) with (\ref{o32wo21-2}) we see that $\omega ^2_{\ 1}\wedge \vartheta ^1=0$, which together with (\ref{o21wth2}) implies $\omega ^2_{\ 1}=0$. We have now $\vartheta ^2=\lambda \, \omega ^3_{\ 1}$ and $\vartheta ^2=\varepsilon _1\, \omega ^3_{\ 1}$, therefore $\lambda =\varepsilon _1=\mbox{const}$.
\bigskip

\section{B\" acklund theorem concerning locally symmetric surfaces}

\begin{theorem}
\label{affBthm}
Let $f:M\to \R ^3$ and $\widehat f:M\to \R ^3$ be
non-degenerate immersions of a two-dimensional connected manifold $M$ into
affine space $\R ^3$, endowed with equiaffine transversal vector fields $\xi $ and $\widehat {\xi }$ respectively.
\smallskip

We denote by $h$ and $\widehat
h$ the corresponding affine fundamental forms, and by $\nu $ and
$\widehat {\nu }$ the conormal maps. 
\medskip

If $f$, $\widehat f$, $\xi $ and $\widehat {\xi }$ satisfy
the following conditions:
\medskip

\noindent \emph{$1^{\circ }$} for every $p\in M$ $f(p)\ne \widehat f(p)$,
the vector $\widehat f(p)-f(p)$ is tangent to $f(M)$ at $f(p)$ and is tangent to $\widehat f(M)$ at $\widehat f(p)$;
\medskip

\noindent \emph{$2^{\circ }$} the spherical representation of $\widehat
f-f$, $M\ni p\mapsto \pi (\widehat f(p)-f(p))\in \Proj ^2(\R )$,
has rank $2$ at every point of $M$;
\medskip

\noindent \emph{$3^{\circ }$} $\det (\widehat f-f,\xi ,\widehat {\xi })\ne 0$ everywhere,
\medskip

\noindent \emph{$4^{\circ }$} the functions $\nu (\widehat {\xi })$ and $\widehat {\nu }(\xi )$ are constant and 
$\nu (\widehat {\xi })\ne 0$ or $\widehat {\nu }(\xi )\ne 0$;
\medskip

\noindent \emph{$5^{\circ }$} $\bigl ( \det ( \widehat f-f,\xi ,\widehat {\xi })\bigr ) ^4\, \cdot \det _{\theta }h\, \cdot \det _{\widehat {\theta }}\widehat h=\bigl ( 1-\nu (\widehat {\xi })\, \widehat {\nu }(\xi )\bigr ) ^4$;
\medskip

\noindent \emph{$6^{\circ }$} for every $Y\in TM$ $\det (f_*(Y),\xi ,\widehat {\xi })=\det (\widehat f_*(Y),\xi ,\widehat {\xi })$;
\medskip

and 
\smallskip

\noindent \emph{$7^{\circ }$} $d(\det (\widehat f-f,\xi,\widehat {\xi }))\wedge d(\det _{\theta }h)= 0$,
\medskip

then affine fundamental forms $h$ and $\widehat h$ are conformal to each other, the connections $\nabla $ and $\widehat {\nabla }$ induced by $(f,\xi )$ and $(\widehat f,\widehat {\xi })$ respectively, are locally symmetric and $\dim \im R=\dim \im \widehat R$.
\end{theorem}
\bigskip

\emph{Proof}. We continue analysis from the proof of Theorem
\ref{wkwhconfdaszh} with the same notation. From $3^{\circ }$ and $5^{\circ }$ we conclude that $1-A\widehat A\ne 0$ and 
\[ \psi (f,\widehat f)=\left ( \frac{1-A\widehat A}{W}\right )
^4\, \frac{1}{H\, \widehat H}=1,\] hence
$h$ and $\widehat h$ are conformal to each other, by Theorem
\ref{wkwhconfdaszh}.
\medskip

It remains to prove that $\nabla R=0$ and $\widehat {\nabla
}\widehat R=0$.
\bigskip

From (\ref{zmformkoneks}) with constant $A$, $\widehat A$
we obtain in particular
\begin{align}
\widehat {\vartheta }^1&=\vartheta ^1+\omega ^1_{\ 1}, \label{dasztheta1} \\
\widehat {\omega }^3_{\ 2}&=-\, \frac{\widehat A(1-A\widehat A)}{W}\, \omega ^2_{\ 2}-\frac{(1-A\widehat A)^2}{W^2}\, \omega ^2_{\ 3}+(\widehat A)^2\, \omega ^3_{\ 2} \nonumber \\
& \ \ \ +\frac{\widehat A(1-A\widehat A)}{W}\, \omega ^3_{\ 3}-\frac{\widehat A(1-A\widehat A)}{W^2}\, dW,
\label{daszomega32} \\
\widehat {\omega }^3_{\ 3}&=(1-A\widehat A)\, \omega ^2_{\
2}-\widehat A\, W\, \omega ^3_{\ 2}-\frac{A(1-A\widehat A)}{W}\,
\omega ^2_{\ 3}+A\widehat A\, \omega ^3_{\ 3}+\frac{1-A\widehat A}{W}\, dW. \label{daszomega33}
\end{align}
For equiaffine vector fields $\xi $ and $\widehat {\xi }$ we have
$\omega ^3_{\ 3}=0$ and $\widehat {\omega }^3_{\ 3}=0$, therefore (\ref{daszomega33}) yields
\begin{equation}
\omega ^2_{\ 2}=\frac{\widehat AW}{1-A\widehat A}\, \omega
^3_{\ 2}+\frac{A}{W}\, \omega ^2_{\ 3}-\frac{dW}{W} \label{omega22}
\end{equation}
and substituting (\ref{omega22}) into (\ref{daszomega32}) we obtain
\begin{equation}
\widehat {\omega }^3_{\ 2}=-\, \frac{1-A\widehat A}{W^2}\, \omega ^2_{\ 3}. \label{daszomega32-2}
\end{equation}
From (\ref{zmformkoneks}) we have also
\begin{align*}
\widehat {\omega }^2_{\ 3}&=-\, AW\, \omega ^2_{\ 2}+A^2\, \omega ^2_{\ 3}-W^2\, \omega ^3_{\ 2}+AW\, \omega ^3_{\ 3}-A\, dW \\ 
&=-\, AW \left ( \omega ^2_{\ 2}-\frac{A}{W}\, \omega ^2_{\ 3}+\frac{dW}{W}\right ) -W^2\, \omega ^3_{\ 2}\\
&=-\, AW\, \frac{\widehat AW}{1-A\widehat A}\, \omega ^3_{\ 2}-W^2\, \omega ^3_{\ 2}
\end{align*}
and it follows that
\begin{equation}
\widehat {\omega }^2_{\ 3}=\frac{-\, W^2}{1-A\widehat A}\, \omega ^3_{\ 2}.
\end{equation}
\medskip

The structural equation (\ref{streqtheta}) with $\vartheta
^3=0$ and $d\vartheta ^3=0$ becomes
\begin{equation}
0=\omega ^3_{\ 1}\wedge \vartheta ^1+\omega ^3_{\ 2}\wedge
\vartheta ^2. \label{streqtheta3}
\end{equation}
Let $\vartheta ^1=s\, \omega ^2_{\ 1}+t\, \omega ^3_{\ 1}$ and
$\omega ^3_{\ 2}=u\, \omega ^2_{\ 1}+v\, \omega ^3_{\ 1}$ with some
functions $s$, $t$, $u$ and $v$. Applying (\ref{streqtheta3}),
(\ref{theta2}) and $\omega ^2_{\ 1}\wedge \omega ^3_{\ 1}\ne 0$
yields
\begin{equation}
s=\frac{\widehat AWu}{1-A\widehat A}+v.
\end{equation}
From (\ref{omega31wedgeomega32}) we obtain 
\[ \omega ^3_{\ 1}\wedge (u\, \omega ^2_{\ 1}+v\, \omega ^3_{\ 1})=H\, (s\, \omega ^2_{\ 1}+t\, \omega ^3_{\ 1})\wedge \left ( -\, \omega ^2_{\ 1}+\frac{\widehat AW}{1-A\widehat A}\, \omega ^3_{\ 1}\right ) ,\]
which implies
\begin{equation}
t=-\, \frac{u}{H}-\frac{\widehat AWs}{1-A\widehat A}=-\, \left ( \frac{1}{H}+\frac{{\widehat A}^2\, W^2}{(1-A\widehat A)^2}\right ) u-\frac{\widehat AW}{1-A\widehat A}\, v.
\end{equation}
Consequently,
\begin{equation}
\vartheta ^1=\left ( \frac{\widehat AWu}{1-A\widehat A}+v\right ) \, \omega ^2_{\ 1}-\left ( \frac{u}{H}+\frac{{\widehat A}^2\, W^2\, u}{(1-A\widehat A)^2}+\frac{\widehat AW\, v}{1-A\widehat A}\right ) \, \omega ^3_{\ 1}. \label{theta1}
\end{equation}
\medskip

We use now the assumption $6^{\circ }$. Since
\begin{align*}
\widehat f_*(Y)-f_*(Y)&=D_Y(\widehat f-f)=D_Yf_*X_1=\omega ^1_{\ 1}(Y)\, f_*X_1+\omega ^2_{\ 1}(Y)\, f_*X_2+\omega ^3_{\ 1}(Y)\, \xi \\
&=\omega ^1_{\ 1}(Y)\, f_*X_1+\omega ^2_{\ 1}(Y)\left ( \frac{A}{W}\, \xi -\frac{1}{W}\, \widehat \xi \right ) +\omega ^3_{\ 1}(Y)\, \xi ,
\end{align*}
\[ \det (\widehat f_*(Y)-f_*(Y),\xi ,\widehat {\xi })=\omega ^1_{\ 1}(Y)\, \det (f_*(X_1),\xi ,\widehat {\xi })=\omega ^1_{\ 1}(Y)\, W,\] 
the equality
\[ \det (f_*(Y),\xi ,\widehat {\xi })-\det (\widehat f_*(Y),\xi ,\widehat {\xi })=0 \]
gives $\omega ^1_{\ 1}=0$ and consequently $\omega ^2_{\ 2}=0$,  because differentiating the equality (\ref{unimod}) we obtain $\omega ^1_{\ 1}+\omega ^2_{\ 2}+\omega ^3_{\ 3}=0$. Similarly from (\ref{unimoddasz}) we obtain 
$\widehat {\omega }^1_{\ 1}+\widehat {\omega }^2_{\ 2}+\widehat {\omega }^3_{\ 3}=0$ and from (\ref{zmformkoneks}) it follows that $\widehat {\omega }^1_{\ 1}=\omega ^1_{\ 1}$, therefore 
$\widehat {\omega }^1_{\ 1}=0$, $\widehat {\omega }^2_{\ 2}=0$ and $\widehat {\vartheta }^1=\vartheta ^1$.
\medskip

Let $\widehat {\omega }^3_{\ 2}=x\, \omega ^2_{\ 1}+y\, \omega ^3_{\ 1}$ with some functions $x$, $y$. Then from
the structural equation $\widehat {\omega }^3_{\ 1}\wedge \widehat {\vartheta }^1+\widehat {\omega }^3_{\ 2}\wedge \widehat {\vartheta }^2=0$, (\ref{daszomega31}), (\ref{theta1}) and (\ref{dasztheta2-2}) we obtain
\[ x=-\, \frac{(1-A\widehat A)^2}{W^2H}\, u.\] 
Using (\ref{daszomega31wedgedaszomega32}) we obtain
\[ y=\frac{1-A\widehat A}{WH}\, \widehat Au-\frac{\widehat HW^3}{(1-A\widehat A)^3}\, \widehat Au-\frac{\widehat HW^2}{(1-A\widehat A)^2}\, v.\]
But $W^4H\widehat H=(1-A\widehat A)^4$, by $4^{\circ }$, hence $y=-\, \frac{(1-A\widehat A)^2}{HW^2}\, v$ and
\begin{equation}
\widehat {\omega }^3_{\ 2}=-\, \frac{(1-A\widehat A)^2}{HW^2}\, \omega ^3_{\ 2}. \label{daszomega32-3}
\end{equation}
Comparing (\ref{daszomega32-3}) with (\ref{daszomega32-2}) we obtain 
\begin{equation}
\omega ^2_{\ 3}=\frac{1-A\widehat A}{H}\, \omega ^3_{\ 2} \label{omega23-vs-omega32}
\end{equation}
and from (\ref{omega22}) with $\omega ^2_{\ 2}=0$
\begin{equation}
dW=\left ( \frac{\widehat AW^2}{1-A\widehat A}+\frac{A(1-A\widehat A)}{H}\right ) \omega ^3_{\ 2}. \label{dW}
\end{equation}
It follows that $\omega ^2_{\ 3}\wedge \omega ^3_{\ 2}=0$. From the fundamental equation
\[ 0=d\omega ^2_{\ 2}=-\, \omega ^2_{\ 1}\wedge \omega ^1_{\ 2}-\omega ^2_{\ 2}\wedge \omega ^2_{\ 2}-\omega ^2_{\ 3}\wedge \omega ^3_{\ 2} \]
we obtain $\omega ^2_{\ 1}\wedge \omega ^1_{\ 2}=0$, which means that 
\begin{equation}
\omega ^1_{\ 2}=\alpha \, \omega ^2_{\ 1} \label{omega12-2}
\end{equation}
for some function $\alpha $. Similarly $\omega ^1_{\ 3}=\beta \, \omega ^3_{\ 1}$, which follows from 
\[ 0=d\omega ^3_{\ 3}=-\, \omega ^3_{\ 1}\wedge \omega ^1_{\ 3}-\omega ^3_{\ 2}\wedge \omega ^2_{\ 3}-\omega ^3_{\ 3}\wedge \omega ^3_{\ 3} .\]
In the same way we obtain $\widehat {\omega }^2_{\ 1}\wedge \widehat {\omega }^1_{\ 2}=0$. From (\ref{zmformkoneks}) we have
\begin{align}
\widehat {\omega }^1_{\ 2}&=\widehat A\, \omega ^1_{\ 2}+\frac{1-A\widehat A}{W}\, \omega ^1_{\ 3},\label{daszomega12} \\
\widehat {\omega }^1_{\ 3}&=-\, W\, \omega ^1_{\ 2}+A\, \omega ^1_{\ 3}. \label{daszomega13}
\end{align}
Using (\ref{daszomega21}) and (\ref{daszomega12}) we obtain
\begin{align}
\widehat {\omega }^2_{\ 1}\wedge \widehat {\omega }^1_{\ 2}&=(A\, \omega ^2_{\ 1}+W\, \omega ^3_{\ 1})\wedge \left ( \widehat A\, \alpha \omega ^2_{\ 1}+\frac{1-A\widehat A}{W}\, \beta \omega ^3_{\ 1} \right ) \nonumber \\
&=\left ( \frac{A(1-A\widehat A)}{W}\, \beta -\widehat AW\, \alpha \right ) \omega ^2_{\ 1}\wedge \omega ^3_{\ 1} \label{powrot}
\end{align}
\smallskip

At first we consider the case $A\ne 0$. It follows that 
\begin{equation}
\beta =\frac{\widehat AW^2}{A(1-A\widehat A)}\, \alpha  \label{beta}
\end{equation}
and
\begin{equation}
\omega ^1_{\ 3}=\frac{\alpha \widehat AW^2}{A(1-A\widehat A)}\, \omega ^3_{\ 1}. \label{omega13-2}
\end{equation}

We have now
\begin{equation}
\widehat {\omega }^1_{\ 2}=\frac{\widehat A}{A}\, \alpha \, \widehat {\omega }^2_{\ 1} \label{daszomega12-2}
\end{equation}
and
\begin{equation}
\widehat {\omega }^1_{\ 3}=\frac{\alpha W^2}{1-A\widehat A}\, \widehat {\omega }^3_{\ 1}. \label{daszomega13-2}
\end{equation}
\medskip

We can already find the curvature tensors of $\nabla $ and $\widehat {\nabla }$. 
\smallskip

We get
\begin{equation}
\nabla _YX_1=\omega ^1_{\ 1}(Y)\, X_1+\omega ^2_{ \ 1}(Y)\,
X_2=\omega ^2_{\ 1}(Y)\, X_2 \label{nablaX1}
\end{equation}
and
\begin{equation}
\nabla _YX_2=\omega ^1_{\ 2}(Y)\, X_1+\omega ^2_{ \ 2}(Y)\,
X_2=\alpha \, \omega ^2_{\ 1}(Y)\, X_1. \label{nablaX2}
\end{equation}
The Gauss equation
\[ d\omega ^k_{ \ l}+\omega ^k_{ \ 1}\wedge \omega ^1_{\ l}+\omega
^k_{ \ 2}\wedge \omega ^2_{ \ l}=-\, \omega ^k_{ \ 3}\wedge \omega
^3_{ \ l}, \ \ \ \ k,l\in \{ 1, 2\} \] now leads to
\begin{equation}
R(X,Y)X_1=-\, \omega ^2_{\ 3}\wedge \omega ^3_{\ 1}(X,Y)\, X_2.
\end{equation}
and
\begin{equation}
R(X,Y)X_2=-\, \omega ^1_{\ 3}\wedge \omega ^3_{\ 2}(X,Y)\, X_1.
\end{equation}
In particular
\begin{equation}
R(X_1,X_2)\, X_1=(1-A\widehat A)\, X_2 \label{R121}
\end{equation}
and
\begin{equation}
R(X_1,X_2)\, X_2=-\, \alpha \frac{\widehat AW^2H}{A(1-A\widehat A)}\, X_1. \label{R122}
\end{equation}
The Ricci tensor is
\begin{align*}
\Ric (X_1,X_1)&=-\, (1-A\widehat A), \\ 
\Ric (X_1,X_2)&=\Ric (X_2,X_1)=0, \\ 
\Ric (X_2,X_2)&=-\, \frac{\alpha \widehat AW^2H}{A(1-A\widehat A)}.
\end{align*}
Applying (\ref{nablaX1}),
(\ref{nablaX2}), (\ref{R121}) and (\ref{R122}) we obtain
\begin{equation}
(\nabla _YR)(X_1,X_2)X_1= (1-A\widehat A)\alpha \left ( \frac{\widehat AW^2H}{A(1-A\widehat A)^2}+1\right ) \omega
^2_{\ 1}(Y)\, X_1
\end{equation}
and
\begin{align}
&(\nabla _YR)(X_1,X_2)\, X_2  \\
&=-\, Y\left ( \frac{\alpha \widehat AW^2H}{A(1-A\widehat A)}\right ) X_1
-(1-A\widehat A)\alpha  \left ( \frac{\widehat AW^2H}{A(1-A\widehat A)^2}+1\right ) \omega ^2_{\ 1}(Y)\, X_2. \nonumber
\end{align}

For $\widehat {\nabla }$ we obtain
\begin{equation}
\widehat {\nabla }_Y\widehat X_1=\widehat {\omega }^2_{\ 1}(Y)\,
\widehat X_2, \ \ \ \ \ 
\widehat {\nabla }_Y\widehat X_2=\frac{\widehat A\, \alpha}{A}\, \widehat {\omega }^2_{\ 1}(Y)\, \widehat X_1, \label{dasznablaX}
\end{equation}
\begin{align*}
\widehat R(\widehat X_1,\widehat X_2)\widehat X_1&=-\, \widehat {\omega }^2_{\ 3}\wedge \widehat {\omega }^3_{\ 1}
(\widehat X_1,\widehat X_2)\, \widehat X_2
=\frac{-\, HW^4}{(1-A\widehat A)^3}\, \widehat {\omega }^3_{\ 2}\wedge \widehat {\omega }^3_{\ 1}(\widehat X_1,\widehat X_2)\, \widehat X_2 \\
&=\frac{W^4H\widehat H}{(1-A\widehat A)^3}\, \widehat X_2=(1-A\widehat A)\, \widehat X_2, \\
\widehat R(\widehat X_1,\widehat X_2)\widehat X_2&=-\, \widehat {\omega }^1_{\ 3}\wedge \widehat {\omega }^3_{\ 2}
(\widehat X_1,\widehat X_2)\, \widehat X_1=\frac{-\, \alpha W^2}{1-A\widehat A}\, \widehat {\omega }^3_{\ 1}\wedge \widehat {\omega }^3_{\ 2}(\widehat X_1,\widehat X_2)\, \widehat X_1 \\
&=\frac{-\alpha W^2\widehat H}{1-A\widehat A}\, \widehat X_1=-\, \frac{\alpha \, (1-A\widehat A)^3}{W^2H}\, \widehat X_1,
\end{align*}
\begin{align*}
\widehat {\Ric }(\widehat X_1,\widehat X_1)&= -\, (1-A\widehat A), \\
\widehat {\Ric }(\widehat X_1,\widehat X_2)&=\widehat {\Ric }(\widehat X_2,\widehat X_1)=0, \\
\widehat {\Ric }(\widehat X_2,\widehat X_2)&=\frac{-\, \alpha \, (1-A\widehat A)^3}{W^2H},
\end{align*}
\begin{align*}
(\widehat {\nabla }_Y\widehat R)(\widehat X_1,\widehat X_2)\widehat X_1&=\frac{\alpha \, (1-A\widehat A)^3}{W^2H}\left ( \frac{\widehat AW^2H}{A(1-A\widehat A)^2}+1\right ) \widehat {\omega }^2_{\ 1}(Y)\, \widehat X_1, \\
(\widehat {\nabla }_Y\widehat R)(\widehat X_1,\widehat X_2)\widehat X_2&=-\, (1-A\widehat A)^3Y\left ( \frac{\alpha }{W^2H}\right ) \widehat X_1 \\ 
& \ \ \ -\frac{\alpha \, (1-A\widehat A)^3}{W^2H}\left (  \frac{\widehat AW^2H}{A(1-A\widehat A)^2}+1\right ) \widehat {\omega }^2_{\ 1}(Y)\, \widehat X_2.
\end{align*}
\bigskip

Next we want to use the assumption $7^{\circ }$: $dW\wedge dH=0$. 
\smallskip

Diferentiating (\ref{omega23-vs-omega32}) we obtain
\[ d\omega ^2_{\ 3}=-\, \frac{1-A\widehat A}{H^2}\, dH\wedge \omega ^3_{\ 2}+\frac{1-A\widehat A}{H}\, d\omega ^3_{\ 2}.\] 
From the fundamental equations and from (\ref{omega12-2}) and (\ref{omega13-2}) we get
\begin{align*}
d\omega ^2_{\ 3}&=-\, \omega ^2_{\ 1}\wedge \omega ^1_{\ 3}= -\, \frac{\alpha \widehat AW^2}{A(1-A\widehat A)}\, \omega ^2_{\ 1}\wedge \omega ^3_{\ 1}, \\
d\omega ^3_{\ 2}&=-\, \omega ^3_{\ 1}\wedge \omega ^1_{\ 2}=\alpha \, \omega ^2_{\ 1}\wedge \omega ^3_{\ 1}.
\end{align*}
It follows that
\begin{equation}
\frac{dH}{H}\wedge \omega ^3_{\ 2}=\alpha \left ( \frac{\widehat AW^2H}{A(1-A\widehat A)^2}+1\right ) \omega ^2_{\ 1}\wedge \omega ^3_{\ 1} 
\end{equation}
and consequently, by (\ref{dW}),
\begin{equation}
dH\wedge dW=A(1-A\widehat A)\, \alpha \left ( \frac{\widehat AW^2H}{A(1-A\widehat A)^2}+1\right ) ^2\omega ^2_{\ 1}\wedge \omega ^3_{\ 1}. \label{dHdW}
\end{equation}
\medskip

If $\widehat A=0$ (and still $A\ne 0$), then (\ref{dHdW}) and $dH\wedge dW=0$ imply $\alpha \equiv 0$.
\smallskip

If $\widehat A\ne 0$ we may compute $d\alpha $ in the following way.

Differentiating (\ref{omega12-2}) and (\ref{omega13-2}) we obtain
\begin{align*}
d\omega ^1_{\ 2}&=d\alpha \wedge \omega ^2_{\ 1}+\alpha \, d\omega ^2_{\ 1}, \\
d\omega ^1_{\ 3}&=\frac{\widehat AW^2}{A(1-A\widehat A)}\, d\alpha\wedge \omega ^3_{\ 1}+\frac{2\alpha \widehat AW}{A(1-A\widehat A)}\, dW\wedge \omega ^3_{\ 1}+\frac{\alpha \widehat AW^2}{A(1-A\widehat A)}\, d\omega ^3_{\ 1}
\end{align*}
and next, after using the fundamental equations, (\ref{dW}) and $\omega ^3_{\ 2}=u\, \omega ^2_{\ 1}+v\, \omega ^3_{\ 1}$,
\begin{align*}
d\alpha \wedge \omega ^2_{\ 1}&=\alpha \left ( \frac{\widehat AW^2H}{A(1-A\widehat A)^2}+1\right ) \frac{1-A\widehat A}{H}\, u\, \omega ^2_{\ 1}\wedge \omega ^3_{\ 1}, \\
d\alpha \wedge \omega ^3_{\ 1}&=-\, \alpha \left ( \frac{\widehat AW^2H}{A(1-A\widehat A)^2}+1\right ) \frac{1-A\widehat A}{H}\left ( \frac{2A}{W}\, u+\frac{A(1-A\widehat A)}{\widehat AW^2}\, v\right ) \omega ^2_{\ 1}\wedge \omega ^3_{\ 1}.
\end{align*}
It follows that
\begin{equation}
d\alpha =-\, \alpha \left ( \frac{\widehat AW^2H}{A(1-A\widehat A)^2}+1\right ) \frac{1-A\widehat A}{H} \left [ \left ( \frac{2A}{W}\, u+\frac{A(1-A\widehat A)}{\widehat AW^2}\, v\right ) \omega ^2_{\ 1}+u\, \omega ^3_{\ 1}\right ] . 
\label{dalpha}
\end{equation}
\bigskip

From (\ref{dHdW}) and $dH\wedge dW=0$ it follows that $\alpha \left ( \frac{\widehat AW^2H}{A(1-A\widehat A)^2}+1\right ) \equiv 0$ on $M$.
Then from (\ref{dalpha}) we conclude that $\alpha $ is constant, because $M$ is connected. 
\bigskip

Now we consider the case $A=0$. Then, by assumption $4^{\circ }$, $\widehat A\ne 0$. We return to (\ref{powrot}) and obtain $\alpha \equiv 0$.
\bigskip

Thus in each case $\alpha =\mbox{const}$.
\bigskip 

If $\alpha =0$, then $\im R_p=\R \, (X_2)_p$, $\im \widehat R_p=\R \, ({\widehat X}_2)_p$, $\dim \im R=\dim \im \widehat R=1$ and $\sign \, \Ric =\sign \, \widehat {\Ric }=-\, \sign (1-A\widehat A)$.
\smallskip

Let $\alpha \ne 0$. Then $\frac{\widehat AW^2H}{A(1-A\widehat A)^2}+1\equiv 0$, which implies 
\[ H=-\, \frac{A(1-A\widehat A)^2}{\widehat AW^2}.\]  
From (\ref{dW}) it follows that $W$ is constant. This clearly forces $H$ to be constant. 
\smallskip

In both cases ($\alpha =0$, $\alpha \ne 0$) we obtain $\nabla R=0$ and $\widehat {\nabla }\widehat R=0$.
\bigskip

We shall show that the case of $\alpha \ne 0$ corresponds to the situation described in the classical B\" acklund theorem or in the B\" acklund theorem for surfaces in Minkowski space. 

\begin{theorem}
If $f$, $\widehat f$, $\xi $, $\widehat {\xi }$ satisfy the assumptions of Theorem \emph{\ref{affBthm}} and the induced connections $\nabla $, $\widehat {\nabla }$ satisfy the condition $\dim \im R=\dim \im \widehat R=2$, then 
$\det (\widehat f-f,\xi, \widehat {\xi })$ and $\det _{\theta }h$ are constant, 
$\R \xi $ and $\R \widehat {\xi }$ are the corresponding affine normals 
and there exists a scalar or pseudoscalar product on $\R ^3$ such that $\xi $ and $\widehat {\xi }$ are orthogonal to the corresponding surfaces with constant, non-zero, length. Moreover the length of $\widehat f-f$ is constant,  
the angle between $\xi $ and $\widehat {\xi }$ is constant and $f$ and $\widehat f$ have the same constant sectional curvature. 
\end{theorem}

\emph{Proof.} We define $G_p\in (\R ^3)^*$ by the equalities
\begin{align*}
G_p(f_*(X_1)_p,f_*(X_1)_p):&=-\, \delta \, (1-A\widehat A), \\
G_p(f_*(X_1)_p,f_*(X_2)_p):&=0, \\
G_p(f_*(X_2)_p,f_*(X_2)_p):&=\delta \, \alpha \, (1-A\widehat A), \\
G_p(f_*X_p,\xi _p):&=0, \\
G_p(\xi _p,\xi _p):&=\delta \, \alpha \, \frac{\widehat A}{A}\, W^2
\end{align*}
with some $\delta \in \{ 1, -1\}$.
We have 
\begin{align*}
SY&=-\, \omega ^1_{\ 3}(Y)\, X_1-\omega ^2_{\ 3}(Y)\, X_2 \\
&=-\, \frac{\alpha\, \widehat AW^2}{A(1-A\widehat A)}\, \omega ^3_{\ 1}(Y)X_1-\frac{1-A\widehat A}{H}\, \omega ^3_{\ 2}(Y)\, X_2 \\
&=-\, \frac{\alpha\, \widehat AW^2}{A(1-A\widehat A)}\, \omega ^3_{\ 1}(Y)X_1+\frac{\widehat AW^2}{A(1-A\widehat A)}\, \omega ^3_{\ 2}(Y)\, X_2 \\
&=\frac{\widehat AW^2}{A(1-A\widehat A)}\left ( -\, \alpha \, \omega ^3_{\ 1}(Y)\, X_1+\omega ^3_{\ 2}(Y)\, X_2\right )
\end{align*}
and
\begin{align*}
&G(f_*X,f_*SY)= \\
&\frac{\widehat AW^2}{A(1-A\widehat A)}\, G(\omega ^1(X)\, f_*(X_1)+\omega ^2(X)\, f_*(X_2),-\alpha \, \omega ^3_{\ 1}(Y)f_*(X_1)+\omega ^3_{\ 2}(Y)\, f_*(X_2))\\
&=\frac{\widehat AW^2}{A(1-A\widehat A)}\left ( \omega ^1(X)\omega ^3_{\ 1}(Y)+\omega ^2(X)\omega ^3_{\ 2}(Y)\right ) \alpha \, \delta \, (1-A\widehat A)
=\delta \,\alpha \, \frac{\widehat A}{A}\, W^2\, h(Y,X).
\end{align*}
Now it is easy to check that $DG=0$, hence we have well defined scalar product on $\R ^3$, which also will be denoted by $G$. The Riemannian or pseudo-Riemannian metric $g$ induced on $M$ by $f$, $g(X,Y)=G(f_*(X),f_*(Y))$, has the sectional curvature
\begin{align*}
\kappa &=\frac{g(R(X_1,X_2)X_2,X_1)}{g(X_1,X_1)\, g(X_2,X_2)-g(X_1,X_2)\, g(X_1,X_2)} \\
&=\frac{g(\alpha \, (1-A\widehat A)\, X_1,X_1)}{-\, \alpha \, \delta ^2\, (1-A\widehat A)^2}=\frac{-\, \delta \, \alpha \, (1-A\widehat A)^2}{-\, \alpha \, \delta ^2\, (1-A\widehat A)^2}=\delta 
\end{align*}
and the same curvature has the metric $\widehat g$ induced by $\widehat f$
\begin{align*}
\widehat {\kappa }&=\frac{\widehat g(\widehat R(\widehat X_1,\widehat X_2)\widehat X_2,\widehat X_1)}{\widehat g(\widehat X_1,\widehat X_1)\, \widehat g(\widehat X_2,\widehat X_2)-\widehat g(\widehat X_1,\widehat X_2)\, \widehat g(\widehat X_1,\widehat X_2)} \\
&=\frac{\widehat g\left ( -\, \frac{\alpha \, (1-A\widehat A)^3}{W^2H}\, \widehat X_1,\widehat X_1\right ) }{-\, \alpha \, \delta ^2\, (1-A\widehat A)^2\frac{\widehat A}{A}}=\frac{A(1-A\widehat A)}{\widehat AW^2H}\, \widehat g(\widehat X_1,\widehat X_1) \\
&=-\, \frac{A(1-A\widehat A)}{\widehat AW^2H}\, \delta \, (1-A\widehat A)=\delta
\end{align*}
because
\begin{align*}
\widehat g(\widehat X_1,\widehat X_1)&=G(\widehat f_*(\widehat X_1),\widehat f_*(\widehat X_1))=G(f_*(X_1),f_*(X_1))=-\, \delta \, (1-A\widehat A), \\
\widehat g(\widehat X_1,\widehat X_2)&=G(\widehat f_*(\widehat X_1),\widehat f_*(\widehat X_2))=
G\left ( f_*(X_1),\widehat A\, f_*(X_2)+\frac{1-A\widehat A}{W}\, \xi \right ) =0, \\
\widehat g(\widehat X_2,\widehat X_2)&=G(\widehat f_*(\widehat X_2),\widehat f_*(\widehat X_2)) \\
&=G\left ( \widehat A\, f_*(X_2)+\frac{1-A\widehat A}{W}\, \xi ,\widehat A\, f_*(X_2)+\frac{1-A\widehat A}{W}\, \xi \right ) \\
&=\widehat A^2\, \delta \, \alpha \, (1-A\widehat A)+\frac{(1-A\widehat A)^2}{W^2}\, \delta \, \alpha \, \frac{\widehat A}{A}\, W^2=\delta \, \alpha \, (1-A\widehat A)\, \frac{\widehat A}{A}
\end{align*}
and
$\frac{A(1-A\widehat A)^2}{\widehat AW^2H}=-\, 1$.
\medskip

We compute
\begin{align*}
G(\widehat f-f,\widehat f-f)&=G(f_*(X_1),f_*(X_1))=-\, \delta (1-A\widehat A) \\
G(\xi , \widehat {\xi })&=G(\xi ,-\, Wf_*(X_2)+A\, \xi )=A\, G(\xi ,\xi )=\delta \, \alpha \, \widehat A\, W^2, \\
G(\widehat {\xi },\widehat {\xi })&=G(-\, Wf_*(X_2)+A\, \xi ,-\, Wf_*(X_2)+A\, \xi ) \\
&=W^2\, \delta \, \alpha \, (1-A\widehat A)+A^2\, \delta \, \alpha \, \frac{\widehat A}{A}\, W^2=\delta \, \alpha \, W^2.
\end{align*}
\medskip

There are five possibilities and we will consider the corresponding cases separately.
\smallskip

({\it i}) \emph{Euclidean case}
\smallskip

If $0<A\widehat A<1$ and $\alpha <0$, then we take $\delta =-\, 1$ and obtain positively definite $G$. Then the square of the length $L$ of $\widehat f-f$ is equal to the positive constant $1-A\widehat A$ and the angle $\measuredangle (\xi ,\widehat {\xi })$ between $\xi $ and $\widehat {\xi }$ is constant too, with
\[ \cos \measuredangle (\xi ,\widehat {\xi })=\sign \widehat A\cdot  \sqrt{A\widehat A}. \] 
Note that 
\begin{align*}
-\, \frac{\sin ^2(\measuredangle (\xi ,\widehat {\xi }))}{L^2}
=-\, \frac{1-\cos ^2(\measuredangle (\xi ,\widehat {\xi }))}{L^2}
=-\, \frac{1-A\widehat A}{1-A\widehat A}=-\, 1=\delta =\kappa =\widehat {\kappa }. 
\end{align*}
\medskip

({\it ii}) \emph{Lorentzian case with timelike congruence $\widehat f-f$ and timelike focal surfaces $f$ and $\widehat f$}
\smallskip

If $0<A\widehat A<1$ and $\alpha >0$, then we take $\delta =1$. We obtain $G(\widehat f-f,\widehat f-f)=-\, (1-A\widehat A)=:-\, L^2$. The plane spanned by $\xi _p$ and $\widehat {\xi }_p$ is spacelike, hence
\[ \cos \measuredangle (\xi ,\widehat {\xi })=\frac{G(\xi ,\widehat {\xi })}{\sqrt{G(\xi ,\xi )}\sqrt{G(\widehat {\xi },\widehat {\xi})}} 
=\sign \widehat A\cdot  \sqrt{A\widehat A}. \]
We obtain 
\[ \frac{\sin ^2(\measuredangle (\xi ,\widehat {\xi }))}{L^2}=1=\delta =\kappa =\widehat {\kappa }.\] 
This case corresponds to (A) of Theorem 2.2 in \cite{AB}.
\medskip

({\it iii}) \emph{Lorentzian case with spacelike congruence $\widehat f-f$ and timelike focal surfaces $f$ and $\widehat f$}
\smallskip

If $A\widehat A>1$ and $\alpha >0$, then we take $\delta =1$ and obtain $G(\widehat f-f,\widehat f-f)=-\, (1-A\widehat A)=:L^2$. Both $\xi _p$ and $\widehat {\xi }_p$ are spacelike, but the plane $\linsp \{ \xi _p, \widehat {\xi }_p\} =\linsp \{ f_*({X_2}_p),\xi _p\} $ is timelike. The hyperbolic angle $\measuredangle (\xi ,\widehat {\xi })$ between two spacelike vectors satisfies the equality 
\[ \cosh ^2(\measuredangle (\xi ,\widehat {\xi }))=\frac{\left ( G(\xi ,\widehat {\xi })\right ) ^2}{G(\xi ,\xi )G(\widehat {\xi },\widehat {\xi })},\] 
which follows from the definition given in \cite{NP-TV}. We obtain $\cosh ^2(\measuredangle (\xi ,\widehat {\xi }))=A\widehat A$ and
\[ \frac{\sinh ^2(\measuredangle (\xi ,\widehat {\xi }))}{L^2}=\frac{\cosh ^2(\measuredangle (\xi ,\widehat {\xi }))-1}{L^2}=\frac{A\widehat A-1}{-(1-A\widehat A)}=1=\delta =\kappa =\widehat {\kappa }.\]
This case corresponds to (B) of Theorem 2.2 in \cite{AB}.
\medskip

({\it iv}) \emph{Lorentzian case with spacelike congruence $\widehat f-f$ and spacelike focal surfaces $f$ and $\widehat f$}
\smallskip

If $A\widehat A>1$ and $\alpha <0$, then we take $\delta =1$. We have $G(\widehat f-f,\widehat f-f)=-\, (1-A\widehat A)=:L^2$ as before, the hyperbolic angle between two timelike vectors satisfies the same equality as above and we obtain 
again
\[ \frac{\sinh ^2(\measuredangle (\xi ,\widehat {\xi }))}{L^2}=1=\delta =\kappa =\widehat {\kappa }.\]
This result is in contradiction with that of Theorem 2.1 in \cite{AB}, where the curvature was claimed to be negative. 
(It seems that in \cite{AB} there is a mistake in going from (2.18) to (2.19), probably $d\omega _{13}$ and $d\omega _{23}$ were incorrect. Moreover, (2.9) on page 43 is in contradiction with $K=-\, \det h_{ij}$ on page 44.)
\medskip

({\it v}) \emph{Lorentzian case with spacelike congruence $\widehat f-f$ and focal surfaces $f$ and $\widehat f$ of different kinds}
\smallskip

If $A\widehat A<0$, then we take $\delta =-\, 1$. Now $G(\widehat f-f,\widehat f-f)=1-A\widehat A=:L^2$ is positive, whereas $G(\xi ,\xi )$ and $G(\widehat {\xi },\widehat {\xi })$ have opposite signs, because $\frac{\widehat A}{A}<0$. 
According to the definition of the hyperbolic angle between timelike vector and spacelike vector, given in \cite{NP-TV}, 
$\measuredangle (\xi ,\widehat {\xi })$ satisfies now the equality
\[ \sinh ^2(\measuredangle (\xi ,\widehat {\xi }))=\frac{\left ( G(\xi ,\widehat {\xi })\right ) ^2}{-\, G(\xi ,\xi )G(\widehat {\xi },\widehat {\xi })}.\] 
We obtain $\sinh ^2(\measuredangle (\xi ,\widehat {\xi }))=-\, A\widehat A$ and
\[ -\, \frac{\cosh ^2(\measuredangle (\xi ,\widehat {\xi }))}{L^2}=-\, \frac{1+\sinh ^2(\measuredangle (\xi ,\widehat {\xi }))}{L^2}=-\, \frac{1-A\widehat A}{1-A\widehat A}=-\, 1=\delta =\kappa =\widehat {\kappa }.\]
The B\" acklund theorem for surfaces of different kinds in Minkowski space can be found in \cite{MR}
\bigskip

{\bf Remark.} In case when both $\nu (\widehat {\xi })$ and $\widehat {\nu }(\xi )$ both equal zero we obtain $W=\mbox{const}$, $dW=0$, the assumption $7^{\circ }$ is satisfied, but we get therefrom no information about relation between $\alpha $ and $\beta $. This case may be characterized by the following proposition.
\medskip

\begin{proposition}
\label{A00}
Let $f$, $\widehat f$, $\xi $, $\widehat {\xi }$ satisfy assumptions $1^{\circ }$, $2^{\circ }$, $3^{\circ }$, $5^{\circ }$ and $6^{\circ }$ of Theorem \emph{\ref{affBthm}} and let $\nu (\widehat {\xi })\equiv 0$ and $\widehat {\nu }(\xi )\equiv 0$. Then there exist local coordinates $x,y$ and functions $H=\det _{\theta }h$, $\alpha $, $\beta $, $\gamma $ satisfying the system of equations
\begin{align}
\alpha &=W^2\, H_y\, e^{-2\gamma }\, \gamma _y+W^2H(e^{-2\gamma }\, \gamma _y)_y+(e^{2\gamma }\, \gamma _x)_x, 
\label{alpha00} \\
\beta &=-\, W^2(e^{-2\gamma }\, \gamma _y)_y-\frac{1}{H}\, (e^{2\gamma }\, \gamma _x)_x+\frac{H_x}{H^2}\, e^{2\gamma }\, \gamma _x,  \label{beta00} \\
\alpha _y&=(\alpha +\beta H)\, \gamma _y, \\
\beta _x&=-\, \frac{1}{H}\, (\alpha +\beta H)\, \gamma _x, 
\end{align}
such that $\vartheta ^i$, $\widehat {\vartheta }^i$, $\omega ^j_{\ k}$ and $\widehat {\omega }^j_{\ k}$ have the following form
\begin{align*}
&\vartheta ^1=\widehat {\vartheta }^1=d\gamma =\gamma _x\, dx+\gamma _y\, dy, \\
&\vartheta ^2=e^{-\, \gamma }\, dx, \ \ \ \widehat {\vartheta }^2=e^{\gamma }\, dy, \\
&\omega ^2_{\ 1}=-\, e^{-\, \gamma }\, dx, \ \ \ \omega ^3_{\ 1}=\frac{e^{\gamma }}{W}\, dy, \ \ \ 
\omega ^1_{\ 2}=-\, \alpha \, e^{-\, \gamma}\, dx, \ \ \ \omega ^1_{\ 3}=\frac{\beta \, e^{\gamma }}{W}\, dy, \\
&\omega ^3_{\ 2}=H\, \omega ^2_{\ 3}=HWe^{-\, 2\gamma }\, \gamma _y\, dx-\frac{e^{2\gamma }}{W}\, \gamma _x\, dy, \\
&\widehat {\omega }^2_{\ 1}=e^{\gamma }\, dy, \ \ \ \widehat {\omega }^3_{\ 1}=\frac{e^{-\, \gamma }}{W}\, dx, \ \ \ 
\widehat {\omega }^1_{\ 2}=\frac{\beta \,e^{\gamma}}{W^2}\, dy, \ \ \ \widehat {\omega }^1_{\ 3}=\alpha \, W e^{-\, \gamma }\, dx, \\
&\widehat {\omega }^2_{\ 3}=HW^4\, \widehat {\omega }^3_{\ 2}=-\, HW^3\, e^{-\, 2\gamma}\, \gamma _y\, dx+W\, e^{2\gamma }\, \gamma _x\, dy. 
\end{align*}
Moreover, $\gamma _x\ne 0$, $\gamma _y\ne 0$ and $W=\det (\widehat f-f,\xi ,\widehat {\xi })$ is a non-zero constant.
\smallskip

The connection $\nabla $ is locally symmetric if and only in $\alpha $ is constant, and $\widehat {\nabla }$ is locally symmetric if and only if $\beta $ is constant. 
\end{proposition} 

Note that from (\ref{alpha00}) and (\ref{beta00}) we obtain
\begin{equation}
\alpha +\beta H=W^2\, H_y\, e^{-2\gamma }\, \gamma _y+\frac{H_x}{H}\, e^{2\gamma }\, \gamma _x.
\label{alfa-plus-betaH}
\end{equation} 
\medskip

\emph{Proof.} If we insert $A=0$ and $\widehat A=0$ into (\ref{dasztheta3-2}) -- (\ref{theta2-3}) and 
(\ref{daszomega32-2}) -- (\ref{daszomega13}), then we obtain 
\begin{align*}
&\vartheta ^2+\omega ^2_{\ 1}=0, \ \ \ \widehat {\vartheta }^2=W\, \omega ^3_{\ 1}=\widehat {\omega }^2_{\ 1}, \ \ \ 
\vartheta ^2=W\, \widehat {\omega }^3_{\ 1}, \\
&\omega ^1_{\ 2}=\alpha \, \omega ^2_{\ 1}, \ \ \ \omega ^1_{\ 3}=\beta \, \omega ^3_{\ 1}, \ \ \ 
\widehat {\omega }^1_{\ 2}=\frac{1}{W}\, \omega ^1_{\ 3}, \ \ \ \widehat {\omega }^1_{\ 3}=-\, W\, \omega ^1_{\ 2}, \\
&\omega ^3_{\ 2}=u\, \omega ^2_{\ 1}+v\, \omega ^3_{\ 1}, \ \ \ \widehat {\omega }^2_{\ 3}=-\, W^2\, \omega ^3_{\ 2}, \ \ \ \widehat {\omega }^3_{\ 2}=\frac{-\, 1}{HW^2}\, \omega ^3_{\ 2}, \ \ \ 
\omega ^2_{\ 3}=\frac{1}{H}\, \omega ^3_{\ 2}, \\
&\widehat {\vartheta }^1=\vartheta ^1=v\, \omega ^2_{\ 1}-\frac{u}{H}\, \omega ^3_{\ 1}, \ \ \ \omega ^1_{\ 1}=\omega ^2_{\ 2}=\widehat {\omega }^1_{\ 1}=\widehat {\omega }^2_{\ 2}=0, \\
&W^4\, H\, \widehat H=1, \ \ \ dW=0. 
\end{align*}
From structural equations with $\omega ^2_{\ 1}=-\, \vartheta ^2$ and $\omega ^1_{\ 2}=\alpha \, \omega ^2_{\ 1}$ it follows that $d\vartheta ^1=0$. 
Hence locally there exists function $\gamma $ such that $\vartheta ^1=d\gamma $. 
It is easy to check that $d(e^{\gamma }\, \vartheta ^2)=0$ and $d(e^{-\, \gamma }\, \widehat {\vartheta }^2)=0$. Moreover $(e^{\gamma }\, \vartheta ^2)\wedge (e^{-\, \gamma }\, \widehat {\vartheta }^2)\ne 0$. Therefore there exist local coordinates $x$, $y$ such that $e^{\gamma }\, \vartheta ^2=dx$ and $e^{-\, \gamma }\, \widehat {\vartheta }^2=dy$.
Next we find the basic $1$-forms $\omega ^2_{\ 1}=-\, \vartheta ^2=-\, e^{-\, \gamma }\, dx$ and $\omega ^3_{\ 1}=frac{1}{W}\, \widehat {\vartheta }^2=\frac{e^{\gamma }}{W}\, dy$. Looking at $\vartheta ^1$ we may find $u$ and $v$, and the rest of $1$-forms is easy to obtain. The system of differential equations for $\alpha $, $\beta $, $\gamma $ and $H$ we get from the fundamental equations. Since $\vartheta ^1\wedge \vartheta ^2\ne 0$ and $\widehat {\vartheta }^1\wedge \widehat {\vartheta }^2\ne 0$, we have $\gamma _x\ne 0$ and $\gamma _y\ne 0$.
\smallskip

We have also 
\begin{align*}
R(X_1,X_2)X_1&=X_2, \\
R(X_1,X_2)X_2&=-\, \beta \, H\, X_1, \\
(\nabla _YR)(X_1,X_2)X_1&=-\, (\alpha +\beta \, H)\, \vartheta ^2(Y)\, X_1, \\
(\nabla _YR)(X_1,X_2)X_2&=-\, Y(\beta \, H)\, X_1+(\alpha +\beta \, H)\, \vartheta ^2(Y)\, X_2  
\end{align*}
and
\begin{align*}
\widehat R(\widehat X_1,\widehat X_2)\widehat X_1&=\widehat X_2, \\
\widehat R(\widehat X_1,\widehat X_2)\widehat X_2&=\frac{-\, \alpha }{W^2\, H}\, \widehat X_1, \\
(\widehat {\nabla }_Y\widehat R)(\widehat X_1,\widehat X_2)\widehat X_1&=\frac{\alpha +\beta \, H}{W^2\, H}\, \widehat {\vartheta }^2(Y)\, \widehat X_1, \\
(\widehat {\nabla }_Y\widehat R)(\widehat X_1,\widehat X_2)\widehat X_2&=Y\left ( \frac{\alpha }{W^2\, H}\right ) \, \widehat X_1-\frac{\alpha +\beta \, H}{W^2\, H}\, \widehat {\vartheta }^2(Y)\, \widehat X_2.  
\end{align*}
If $\nabla R=0$ then $\alpha +\beta \, H=0$ and $\beta \, H=\mbox{const}$, hence $\alpha =-\, \beta \, H$ is also constant. Conversely, if $\alpha =\mbox{const}$, then $\alpha _y=0$ and from the system of differential equations we obtain $\alpha +\beta \, H=0$, next $\beta \, H=-\, \alpha =\mbox{const}$ and $\nabla R=0$.
\smallskip

If $\widehat {\nabla }\widehat R=0$, then $\alpha +\beta \, H=0$ and $\frac{\alpha }{H}$ is constant, and now $\beta =\frac{-\, \alpha }{H}$. Conversely, if $\beta $ is constant, then from $\beta _x=0$ we obtain $\alpha +\beta \, H=0$, hence $\frac{\alpha }{H}=-\, \beta $ is constant and $\widehat {\nabla }\widehat R=0$.

\section{The particular case when connections are induced by affine normal vector fields}

\begin{theorem}
\label{affBthm-apn}
Let $f:M\to \R ^3$ and $\widehat f:M\to \R ^3$ be
non-degenerate immersions of a two-dimensional real manifold $M$ into
affine space $\R ^3$.
\smallskip

We denote by $\xi $ and $\widehat {\xi }$ the affine normal vector
field for $f$ and $\widehat f$ respectively, by $h$ and $\widehat
h$ the corresponding affine fundamental forms, and by $\nu $ and
$\widehat {\nu }$ the conormal maps. Let $\varepsilon =\sign \det
h_{ij}$ and $\widehat {\varepsilon }=\sign \det \widehat h_{ij}$.
\medskip

Let $f$ and $\widehat f$ satisfy
the following conditions:
\medskip

\noindent \emph{(i)} for every $p\in M$ $f(p)\ne \widehat f(p)$,
moreover the vector $\widehat f(p)-f(p)$ is tangent to $f(M)$ at $f(p)$ and is tangent to $\widehat f(M)$ at $\widehat f(p)$,
\medskip

\noindent \emph{(ii)} the spherical representation of $\widehat
f-f$, $M\ni p\mapsto \pi (\widehat f(p)-f(p))\in \Proj ^2(\R )$,
has rank $2$ at every point of $M$,
\medskip

\noindent \emph{(iii)} the functions $\nu (\widehat {\xi })$ and $\widehat {\nu }(\xi )$ are constant,
\medskip

\noindent \emph{(iv)} $\det ( \widehat f-f,\xi ,\widehat {\xi })$ is a non-zero constant,
\medskip

\noindent \emph{(v)}
$\bigl | \det ( \widehat f-f,\xi ,\widehat {\xi })\bigr | = \bigl | 1-\nu (\widehat {\xi })\, \widehat {\nu }(\xi )\bigr | $,
\medskip

\noindent \emph{(vi)} $\varepsilon =\widehat {\varepsilon }$.
\medskip

Then affine fundamental forms $h$ and $\widehat h$ are conformal to each other.
\bigskip

If moreover
\smallskip

\noindent \emph{(vii)} $\widehat {\nu }(\xi )+\varepsilon \, \nu
(\widehat {\xi })=0$,
\smallskip

then the Blaschke connections $\nabla $ and $\widehat {\nabla }$, of $f$ and $\widehat f$ respectively, are locally symmetric.
\end{theorem}
\bigskip

\emph{Proof.} Without loss of generality we may assume that $\det ( \widehat f-f,\xi ,\widehat {\xi })=1-\nu (\widehat {\xi })\, \widehat {\nu }(\xi )$, because affine normal vector field $\widehat {\xi }$ may be replaced by $-\, \widehat {\xi }$. We retain our previous notation, so we have now $W=1-A\widehat A$. The case $A=\widehat A=0$ is described in Theorem 1.5 of \cite{MR}. We may also use (\ref{alfa-plus-betaH}) with constant $H$ and next use Proposition \ref{A00}.
\smallskip

If $A\ne 0$ or $\widehat A\ne 0$, then $(f,\xi )$ and $(\widehat f,\widehat {\xi })$
satisfy the assumptions $1^{\circ }$ -- $5^{\circ }$ and $7^{\circ }$ of Theorem \ref{affBthm}.
It suffices to check whether they satisfy $6^{\circ }$.
\smallskip

We will show that the assumption $\widehat A+\varepsilon \, A=0$ implies $\omega ^1_{\ 1}=0$, which is equivalent to $6^{\circ }$.
We proceed as in the first part of the proof of Theorem \ref{affBthm} and obtain the formulae corresponding to (\ref{dasztheta1}), (\ref{dasztheta2-2}), (\ref{theta2-3}), (\ref{daszomega32-2}), (\ref{omega22}) and (\ref{theta2}), when $W=1-A\widehat A=\mbox{const}$:
\begin{align}
\widehat {\vartheta }^1&=\vartheta ^1+\omega ^1_{\ 1}, \label{dasztheta1-apn} \\
\widehat {\vartheta }^2&=\omega ^3_{\ 1}, \label{dasztheta2-apn} \\
\widehat {\omega }^3_{\ 1}&=\vartheta ^2, \label{daszomega31-apn} \\
\widehat {\omega }^3_{\ 2}&=\frac{-\, 1}{W}\, \omega ^2_{\ 3}, \label{daszomega32-1-apn} \\
-\, \omega ^1_{\ 1}&=\widehat A\, \omega ^3_{\ 2}+\frac{A}{W}\, \omega ^2_{\ 3}, \label{11-apn} \\
\vartheta ^2&=-\, \omega ^2_{\ 1}+\widehat A\, \omega ^3_{\ 1}. \label{theta2-apn}
\end{align}
If we bring together (\ref{daszomega32-1-apn}) and (\ref{11-apn}), then we obtain
\begin{equation}
\label{daszomega32-2-apn}
\widehat {\omega }^3_{\ 2}=\frac{\widehat A}{A}\, \omega ^3_{\ 2}+\frac{1}{A}\, \omega ^1_{\ 1}.
\end{equation}
Note, that if $\widehat A+\varepsilon \, A=0$ and $(A,\widehat A)\ne (0,0)$, then $A\ne 0$. 
\smallskip

Substituting (\ref{dasztheta1-apn}), (\ref{dasztheta2-apn}), (\ref{daszomega31-apn}) and (\ref{daszomega32-2-apn}) into
\begin{align*}
&\widehat {\omega }^3_{\ 1}\wedge \widehat {\vartheta }^1+\widehat {\omega }^3_{\ 2}\wedge \widehat {\vartheta }^2=0, \\
&\widehat {\omega }^3_{\ 1}\wedge \widehat {\omega }^3_{\ 2}=\widehat {\varepsilon }\, \widehat {\vartheta }^1\wedge \widehat {\vartheta }^2=\varepsilon \, \widehat {\vartheta }^1\wedge \widehat {\vartheta }^2
\end{align*}
we obtain
\begin{align*}
&\omega ^1_{\ 1}\wedge \left ( -\, \vartheta ^2+\frac{1}{A}\, \omega ^3_{\ 1}\right ) - \left ( \vartheta ^1\wedge \vartheta ^2+\frac{\widehat A}{A}\, \omega ^3_{\ 1}\wedge \omega ^3_{\ 2}\right ) =0, \\
&\omega ^1_{\ 1}\wedge \left ( -\, \frac{1}{A}\, \vartheta ^2-\varepsilon \, \omega ^3_{\ 1}\right ) + \left ( \varepsilon \, \omega ^3_{\ 1}\wedge \vartheta ^1-\frac{\widehat A}{A}\, \omega ^3_{\ 2}\wedge \vartheta ^2\right ) =0.
\end{align*}
But $\frac{\widehat A}{A}=-\, \varepsilon $, therefore
\begin{align*}
\vartheta ^1\wedge \vartheta ^2+\frac{\widehat A}{A}\, \omega ^3_{\ 1}\wedge \omega ^3_{\ 2}&=\varepsilon \left ( \varepsilon \, \vartheta ^1\wedge \vartheta ^2-\omega ^3_{\ 1}\wedge \omega ^3_{\ 2}\right ) =0, \\
\varepsilon \, \omega ^3_{\ 1}\wedge \vartheta ^1-\frac{\widehat A}{A}\, \omega ^3_{\ 2}\wedge \vartheta ^2&=\varepsilon \, \left ( \omega ^3_{\ 1}\wedge \vartheta ^1+\omega ^3_{\ 2}\wedge \vartheta ^2\right ) =0
\end{align*}
and consequently 
\begin{align*}
&\omega ^1_{\ 1}\wedge \left ( -\, \vartheta ^2+\frac{1}{A}\, \omega ^3_{\ 1}\right ) =0, \\
&\omega ^1_{\ 1}\wedge \left ( -\, \frac{1}{A}\, \vartheta ^2-\varepsilon \, \omega ^3_{\ 1}\right ) =0.
\end{align*}
The $1$-forms 
\begin{align*}
-\, \vartheta ^2+\frac{1}{A}\, \omega ^3_{\ 1}&=\omega ^2_{\ 1}+\frac{1-A\widehat A}{A}\, \omega ^3_{\ 1}, \\  
-\, \frac{1}{A}\, \vartheta ^2-\varepsilon \, \omega ^3_{\ 1}&=\frac{1}{A}\, \omega ^2_{\ 1}-\frac{\widehat A+\varepsilon \, A}{A}\, \omega ^3_{\ 1}=\frac{1}{A}\, \omega ^2_{\ 1}
\end{align*} 
are linearly independent, hence $\omega ^1_{\ 1}=0$. 
\medskip

It follows that we may now apply Theorem \ref{affBthm}.

\end{document}